\documentclass[12pt,a4paper,reqno,twoside]{amsart}

\usepackage[english]{babel}
\usepackage{stmaryrd}
\usepackage{dsfont}
\usepackage{bbm}
\usepackage[symbol*,ragged]{footmisc}
\usepackage[colorlinks,linkcolor=red,anchorcolor=blue,citecolor=blue,urlcolor=blue]{hyperref}
\usepackage{color,xcolor}

\usepackage{geometry}
\usepackage{amssymb}
\usepackage{amsmath}
\usepackage{mathrsfs}
\usepackage{amsfonts}
\usepackage{epsfig}

\usepackage{amsthm}
\usepackage{amsxtra}
\usepackage{bbding}
\usepackage{epsfig}
\usepackage{graphicx}
\usepackage{latexsym}
\usepackage{mathbbol}
\usepackage{bbold}

\usepackage{pifont}
\usepackage{wasysym}
\usepackage{skull}
\usepackage{float}

\DeclareSymbolFontAlphabet{\mathbb}{AMSb}
\DeclareSymbolFontAlphabet{\mathbbol}{bbold}

\usepackage{amscd}
\usepackage[all]{xy}
\allowdisplaybreaks[4]
\usepackage{setspace}

\theoremstyle{plain}

\newtheorem*{corollary*}{\normalfont\scshape Corollary}

\theoremstyle{remark}
\newtheorem*{remark*}{\normalfont\scshape Remark}

\numberwithin{equation}{section}
\addtocounter{footnote}{1}

\renewcommand{\footnoterule}{
  \kern -3pt
  \hrule width 2.5in height 0.4pt
  \kern 3pt
}

\makeatletter
\@ifundefined{MakeUppercase}{}{}
\makeatother

\begin{document}
	
\title[On the average number of subgroups of the group ${\Bbb Z}_{m}\times {\Bbb Z}_{n}$]{On the average number of subgroups of the group ${\Bbb Z}_{n_1}\times {\Bbb Z}_{n_2}$ }

\author {Wenguang Zhai}

\address{[Wenguang Zhai]  Department of Mathematics, China University of Mining and Technology,          Beijing 100083, People's Republic of China}

\email{zhaiwg@hotmail.com}

\date{}


\subjclass[2000]{11N37, 11N45}
\keywords{subgroup,
Dirichlet series, asymptotic formula, divisor problem, Tauber theorem}

\thanks{This work is  partially supported  by
the National Natural Science Foundations of China(Grant Nos. 12471009,  12301006)
and  partially supported  by  Beijing Natural Science Foundation (Grant No. 1242003).}

\begin{abstract}
 Let ${\mathbb Z}_{n}$ be the additive group of residue classes modulo $n.$ For any positive integers $n_1$ and $n_2$, let $s(n_1,n_2)$ and $c(n_1,n_2)$ denote the number of subgroups and the number of   cyclic subgroups of the   group ${\mathbb  Z}_{n_1}\times {\mathbb  Z}_{n_2}$, respectively. The aim of this paper is to study the asymptotic behavior of the sums
 $\sum_{n_1,n_2\le x}s(n_1,n_2)$ and $\sum_{n_1,n_2\le x}c(n_1,n_2)$.
  Some  sharper asymptotic  results are given for the these sums. Mean values of the error terms are also studied.
\end{abstract}

\maketitle

\addtocounter{footnote}{1}

\tableofcontents
\section{Introduction}

\subsection{On the average number of subgroups of ${\Bbb Z}_{n_1}\times {\Bbb Z}_{n_2}$}

For any fixed integer $k\geq 2,$ let $\tau_k(n)$ denote the number of ways
$n$ can be written as a product of $k$ natural numbers and
let  $D_k(x) $ denote its summatory function up to $x.$ Especially when
$k=2,$ let $\tau(n)=\tau_2(n)$ and $D(x)=D_2(x).$
Dirichlet first proved
\begin{equation}\label{1.1}
D(x)=x\log x+(2\gamma-1)x+\Delta(x),\ \ \ \  \Delta(x)=O(x^{1/2}),
\end{equation}
where $\gamma$ is the  Euler's constant. The well-known Dirichlet divisor problem
is to study properties of $\Delta(x)$, for example, its upper bound, moments, $\Omega$-results, sign changes, etc. The Dirichlet divisor problem is a classical problem in analytic number theory and there are a vast of results on this subject. See, for example, \cite{Cr, HB, Iv3, LT, T, Ts, Ts2, Zh1}.

For any $n\in {\Bbb N}, $  let
 $ {\mathbb Z}_{n}$ denote the additive group of residue classes modulo $n.$
Let $r\geq 1$ be a fixed integer and $n_1, \cdots, n_r\in {\Bbb N}.$
Consider the set $ {\mathbb Z}_{n_1}\times\cdots \times {\mathbb Z}_{n_r}$.
For any $(a_1,\cdots  a_r)\in {\mathbb Z}_{n_1}\times \cdots \times {\mathbb Z}_{n_r},
(b_1,\cdots  b_r) \in  {\mathbb Z}_{n_1}\times \cdots \times {\mathbb Z}_{n_r},$ define
\begin{eqnarray*}
 (a_1,\cdots  a_r) \oplus  (b_1,\cdots  b_r)
 =(a_1+b_1(mod\ n_1), \cdots, a_r +b_r(mod\ n_r))
   \in {\mathbb Z}_{n_1}\times\cdots \times  {\mathbb Z}_{n_r}. \end{eqnarray*}
It is easy to check that the set $ {\mathbb Z}_{n_1}\times\cdots \times  {\mathbb Z}_{n_r}  $ forms a group  under the operation $\oplus$. Let $s(n_1, \cdots, n_r)$ and $c(n_1, \cdots, n_r)$ denote the number of subgroups and  the number of cyclic subgroups of ${\Bbb Z}_{n_1}\times \cdots \times {\Bbb Z}_{n_r},$ respectively.

If  $r=1$, then $s(n)=c(n)=\tau(n).$   So We have
 $$\sum_{n\leq x}s(n)=\sum_{n\leq x}c(n)=\sum_{n\leq x}\tau(n)=D(x).$$
 This means that we can look at the Dirichlet divisor problem from
 the point of view of group theory.

The case $r=2$ is interesting.
 If $gcd(n_1,n_2)=1,$ then
  $ {\mathbb Z}_{n_1}\times {\mathbb Z}_{n_2} $ is cyclic and isomorphic to ${\mathbb Z}_{n_1n_2}.$ If $gcd(n_1,n_2)>1,$ then ${\mathbb Z}_{n_1}\times {\mathbb Z}_{n_2}  $ has rank two.    T\'oth used an  elementary method in \cite{Fu} to prove the following  compact formulas:
for every $n_1,n_2\in {\Bbb N},$
\begin{eqnarray*}
s(n_1,n_2)&&=\sum_{d|n_1,e|n_2}gcd(d,e)\\
&&=\sum_{d|gcd(n_1,n_2)}\varphi(d)\tau(n_1/d)\tau(n_2/d)\nonumber
\end{eqnarray*}
and
\begin{eqnarray*}
c(n_1,n_2)&&=\sum_{d|n_1,e|n_2}\varphi(gcd(d,e))\\
&&=\sum_{d|gcd(n_1,n_2)}(\mu*\varphi)(d)\tau(n_1/d)\tau(n_2/d).\nonumber
\end{eqnarray*}
So both $s(n_1,n_2)$ and $c(n_1,n_2)$ have closed formulas.

If $r\geq 3$,  then   T\'oth \cite{To} proved the formula
$$c(n_1, \cdots, n_r)=\sum_{d_1|n_1,\cdots, d_r|n_r}\frac{\varphi(d_1)\cdots \varphi(d_r)}{\varphi([d_1, \cdots, d_r])}.$$

However, there is no such a closed formula for $s(n_1,\cdots, n_r)\ (r\geq 3)$.

Suppose $x>0$ is a real number. Define
\begin{eqnarray*}
&&S^{(1)}(x):=\sum_{n_1\leq x, n_2\leq x}s(n_1,  n_2),\ \ \
S^{(2)}(x):=\sum_{\stackrel{n_1\leq x, n_2\leq x}{gcd(n_1, n_2)>1}}s(n_1,n_2),\\
&&C^{(1)}(x):=\sum_{n_1\leq x, n_2\leq x}c(n_1,  n_2),\ \ \
C^{(2)}(x):=\sum_{\stackrel{n_1\leq x, n_2\leq x}{gcd(n_1,  n_r)>1}}c(n_1, n_2).
\end{eqnarray*}
Here $S^{(2)}(x)$ and $C^{(2)}(x)$ denote the number of subgroups and cyclic subgroups of the groups ${\Bbb Z}_{n_1} \times{\Bbb Z}_{n_2}$ having rank  two with $n_1,  n_2\leq x.$ The above four counting functions can be viewed as two dimensional generalizations of
$D(x).$

Nowak and T\'oth \cite{NT} first studied the asymptotic behaviour of these counting functions. Let
\begin{equation}
\mathcal{F}[x]:=\{S^{(1)}(x), S^{(2)}(x), C^{(1)}(x), C^{(2)}(x) \}.
\end{equation}
Nowak and T\'oth  proved the following asymptotic formula
\begin{eqnarray}\label{Nowak-Toth}
&& f(x)=x^2\left(\sum_{r=0}^3A_{r,f}\log^r x\right)+O(x^{\frac{1117}{701}+\varepsilon})
\end{eqnarray}
holds for any $f(x)\in \mathcal{F}[x],$
where $A_{r,f}(r=0,1,2,3)$ are explicit constants.

By using the well-known complex integration approach, T\'oth  and Zhai \cite {TZ1} proved that
the error term $O(x^{\frac{1117}{701}+\varepsilon})$ in the asymptotic formula
(\ref{Nowak-Toth}) can be replaced
 by $O(x^{3/2}(\log x)^{6.5})$ for any $f(x)\in \mathcal{F}[x]$.

T\'oth  and Zhai \cite {TZ2} proved that the asymptotic formula
\begin{equation*}
\sum_{n_1, n_2, n_3\leq x}c(n_1, n_2, n_3)=x^3\left(\sum_{j=0}^7c_j\log^j x\right)+O(x^{8/3+\varepsilon})
\end{equation*}
holds, where $c_j(j=0,\cdots, 7)$ are explicit constants.
When $r\geq 4,$ D. Essouabri, C. Salinas Zavala and
 L. T\'{o}th \cite{EST} proved that the asymptotic formula
\begin{equation*}
\sum_{n_1, \cdots, n_r\leq x}c(n_1, \cdots, n_r)=x^rQ_1(\log x)+O(x^{r-\mu_r})
\end{equation*}
holds for some positive constant $\mu_r>0,$ where $Q_1(u)$ is a polynomial in $u$
of degree $2^r-1.$

The aim of this paper is to prove some sharper results for $r=2.$

\subsection{Statements of results}\

For any $f(x)\in \mathcal{F}[x],$ define
\begin{eqnarray*}
&&E_f(x):=f(x)-x^2\left(\sum_{r=0}^3A_{r,f}\log^r x\right).
\end{eqnarray*}

 Our first result is  the following Theorem 1, which further improves
previous upper bound results of $E_f(x)$.

{\bf Theorem 1.}  {\it Suppose $f(x)\in \mathcal{F}[x]$. Then
\begin{eqnarray*}
&&E_f(x)\ll x^{\frac{139}{96}+\varepsilon}.
\end{eqnarray*}}

{\bf Remark 1.} Numerically, we have

$$1117/701=1.5934\cdots, \ \ 3/2=1.5, \ \ 139/96=1.4479\cdots.$$

Now we study mean values of $E_f(x)$ for any $f(x)\in \mathcal{F}[x].$ The first result in this direction is the following

{\bf Theorem 2.} {\it  Suppose $f(x)\in \mathcal{F}[x]$. Then we have
\begin{eqnarray}
&&\int_1^X E_f(x) dx\ll X^{\frac{13}{6}+\varepsilon}.
\end{eqnarray} }

It is natural to ask the question: what is the best possible upper   bound
of $E_f(x)$ for $f(x)\in \mathcal{F}[x]?$ For this question, we propose the following conjecture.

{\bf Conjecture 1.} {\it  Suppose $f(x)\in \mathcal{F}[x]$. Then
\begin{eqnarray*}
&&E_f(x)\ll x^{\frac{4}{3}+\varepsilon}.
\end{eqnarray*}}

The following Theorem 3 is a result about the weighted mean square of $E_f(x)\ (f(x)\in \mathcal{F}[x])$, which  implies that Conjecture 1 is true on average.

{\bf Theorem 3.} {\it  Suppose $f(x)\in \mathcal{F}[x]$. Then we have the asymptotic formula
\begin{eqnarray}
&&\int_1^X E_f^2(x) x^{-\frac{11}{3}}dx=\left(c_f+O\left(\frac{1}{\log\log X}\right)\right)\log^3 X,
\end{eqnarray}
where $c_f>0  $ is a positive constant.
Furthermore, we have
$$c_{S^{(1)}}=c_{S^{(2)}}=\frac{8B_{2,1}(2/3)}{81\pi^2},\ \ \
c_{C^{(1)}}=c_{C^{(2)}}=\frac{8B_{2,2}(2/3)}{81\pi^2},
$$ where $B_{2,j}(\sigma)\ (j=1,2)$ are defined by (\ref{I2j-cons}).
}

From Theorem 3 we immediately get the following corollary.

{\bf Corollary 1.} {\it  Suppose $f(x)\in \mathcal{F}[x]$. Then we have the estimates
\begin{eqnarray}\label{1.6}
&&\int_1^X E_f^2(x)  dx\ll  \frac{X^{11/3}\log^3 X}{\log\log X} \ \ (X\rightarrow \infty)
\end{eqnarray}
and
\begin{eqnarray}\label{1.7}
&&E_f(x)=\Omega(x^{4/3}\log x)\ \ (x\rightarrow \infty).
\end{eqnarray}
}

We propose the following Conjecture 2 about the mean square of $E_f(x).$

{\bf Conjecture 2.} {\it Suppose $f(x)\in \mathcal{F}[x]$. Then the weak asymptotic formula
 \begin{eqnarray}
&&\int_1^X E_f^2(x)  dx\sim c^{(2)}_f X^{11/3}\log^2 X\ ( X\rightarrow \infty)
\end{eqnarray} holds
for some positive constant $c^{(2)}_f>0.$ More precisely, we should have the asymptotic formula
  \begin{eqnarray}\label{con-asy-2}
&&\int_1^X E_f^2(x)  dx=   X^{11/3}\left(c^{(2)}_f\log^2 X+c^{(1)}_f\log  X+c^{(0)}_f\right)+O(X^{11/3-\eta}),
\end{eqnarray}
 where $0<\eta<1/2 , c^{(2)}_f, c^{(1)}_f, c^{(0)}_f $ are real constants.
}

From Theorems 1, 2, 3   we get the following corollary.

{\bf Corollary 2.} {\it  Suppose $f(x)\in \mathcal{F}[x]$. Then we have
the following assertions.

(1)
$ E_f(x)$ changes its sign infinitely often.

(2) Let $N_f(X)$ denote the number of sign changes of $E_f(x)$ in the interval
$[1, X].$  Then $N_f(X)\gg (\log\log X)^2$ when $X$ tends to infinity.

(3) If the asymptotic formula (\ref{con-asy-2}) were true, then
$N_f(X)\gg \min(X^\eta, X^{5/96-2\varepsilon})  $ when $X$ tends to infinity.
}

{\bf Remark 2.} The sums of the form $\sum_{m,n}s(m,n)$ and
$\sum_{m,n}c(m,n)$ are  sums over integer points in subsets of ${\Bbb N}^2$.
A main ingredient of this paper is to transform such sums over ${\Bbb N}^2$ into
 sums over ${\Bbb N}$, namely, we change a two-dimensional problem  to
a one-dimensional problem. The idea was used to treat $S(x)$ in Section 5,
which is a sum important to all theorems of our paper.

The structure of the paper is as follows. In Section 2 we quote some lemmas which are needed for our proof. In Section 3 we study a special three-dimensional divisor problem, which is a generalization of $D_3(x).$
In Section 4 we study a two-dimensional generalization of the Dirichlet divisor problem, which is useful for $S^{(2)}(x)$ and  $C^{(2)}(x).$
The proof of Theorem 1 and the proof of Theorem 2 will be given in
 Section 5 and  Section 6, respectively. In Section 7 we study some integrals involving the Riemann zeta-function in the critical strip $1/2<\sigma<1$, which is important for the proof of Theorem 3. The proofs of Theorem 3 and two corollaries will be given in Section 8.

{\bf Notation.}
For any positive integer $m>0,$ let ${\mathbb  Z}_{m}$ denote the additive group of residue classes modulo $m.$
Throughout this paper, ${\mathbb  N}$ denotes the set of all positive integers, $\varphi$ denotes Euler's totient function,
 $\mu$ denotes the M\"obius function, $\zeta$
denotes the Riemann zeta-function, $\tau_{\ell}(n)$ denotes the number of ways $n$ can be written as a product of $\ell$ positive integers($\tau(n)=\tau_2(n), \ \tau_1(n)\equiv 1$), $f*g$ means the Dirichlet convolution $f*g(n)=\sum_{n=n_1n_2}f(n_1)g(n_2).$   We always use $\varepsilon$ denote
a sufficiently small positive constant. The symbol $n\sim N$ means $N<n\leq 2N$
 and $n\asymp N$ means that there exist two absolute positive constants $0<c_1<c_2$ such that $c_1N\leq n\leq c_2N.$ As usual, the symbols $f=O(g)$ and $f\ll g$ mean that $|f|\leq Cg$ for some positive constant $C.$ The expression $e(t)$ means  $e(t)=e^{2\pi i t}.$

\section{Some preliminary lemmas}

In order to prove our results, we need the following lemmas.

{\bf Lemma 2.1.} {\it The Riemann-zeta function $\zeta(s)$
satisfies the functional equation
\begin{equation}
\zeta(s)=\chi(s)\zeta(1-s),
\end{equation}
where
$$\chi(s)=\frac{\Gamma((1-s)/2)}{\Gamma(s/2)}\pi^{s-\frac 12}.$$
}

\begin{proof}
See, for example,   Chapter 1 of Ivi\'{c} \cite{Iv}, or Remark 1.2 in   \cite {Iv2}.
\end{proof}

{\bf Lemma 2.2.} {\it  Suppose $-10\leq u\leq 10, v\geq 1.$ Then we have
\begin{eqnarray}
&&\chi(u+iv)=\left(\frac{v}{2\pi}\right)^{\frac 12 -u-iv}e^{i(v+\frac \pi 4)}
\left(1+O\left(\frac 1v\right)\right),\\
&&\frac{\chi^{\prime}(u+iv)}{\chi(u+iv)}=-\log\frac{v}{2\pi}+O\left(\frac 1v\right).\nonumber
\end{eqnarray}
}

\begin{proof}
This is formula (2.17) of \cite {Iv2}.
\end{proof}

{\bf Lemma 2.3.} {\it   Suppose $j\geq 0$ is a fixed integer and $ |t|\geq 2$.
Let $\zeta^{(j)}(s)$ denote the $j$-th derivative of $\zeta(s).$ Then one has
 \begin{eqnarray}
\zeta^{(j)}(\sigma+it)\ll \left\{\begin{array}{ll}
|t|^{1/2-\sigma} \log^{1+j} |t|,&\mbox{$-1/2 \leq \sigma\leq 0,$}\\
|t|^{1/2-2\sigma/3}\log^{1+j} |t|,& \mbox{$ 0\leq \sigma\leq 1/2,$}\\
 |t|^{(1-\sigma)/3}\log^{1+j} |t|,&\mbox{$ 1/2\leq \sigma\leq 1 ,$}\\
\min(\log^{1+j} |t|,1/(\sigma-1)^{1+j}),&\mbox{$\sigma> 1.$}
\end{array}\right.
\end{eqnarray}
}

\begin{proof}
The estimate for the case $j=0$ follows from the classical bounds
$$\zeta(it)\ll  (|t| + 2)^{1/2}\log(|t|+2), \ \ \zeta(1/2 + it)\ll
 (|t| + 2)^{1/6}, \ \   \zeta(1 + it)\ll  \log(|t| + 2)$$
and the Phragmen-Lindel\"{o}f principle. The estimate for the case $j\geq 1$ follows from the
result of the case $j=0$ and Cauchy's theorem.
\end{proof}

{\bf Lemma 2.4.} {\it   Let $j\geq 0$ be a fixed integer and
  $\zeta^{(j)}(s)$ denote the $j$-th derivative of $\zeta(s).$  Then we have the estimate
\begin{eqnarray*}
&&\int_1^T |\zeta^{(j)}(7/12+it)|^6\ll T^{1+\varepsilon}.
\end{eqnarray*}
  }

\begin{proof}
The case $j=0$ is contained in Theorem 8.4 of Ivi\'c \cite{Iv}.
The case $j\geq 1$ follows from the result of the case $j=0$ and  Cauchy's Theorem.
\end{proof}

{\bf Lemma 2.5.}  {\it Suppose that $f(x)$ and $g(x)$ are real-valued functions on the interval $[a,b]$
which satisfy the conditions

1) $f^{(4)}(x)$ and $g^{\prime\prime}(x)$ are continuous;

2) there exist  positive numbers $H, U, A, $ $  0<b-a\leq U,$ such that
$$A^{-1}\ll f^{\prime\prime}(x)\ll A^{-1},\ \ f^{(3)}(x)\ll A^{-1}U^{-1}, \ \ f^{(4)}(x)\ll A^{-1}U^{-2},$$

$$g(x)\ll H, \ \ g^{\prime}(x)\ll HU^{-1},\ \ g^{\prime\prime}(x)\ll HU^{-2};$$

3) $f^{\prime}(c)=0$ for some $c$, $a\leq c\leq b.$

Then
\begin{eqnarray*}
&&\int_a^b g(x)e(f(x))dx=\frac{1+i}{\sqrt 2}\cdot \frac{g(c)e(f(c))}{\sqrt {f^{\prime\prime}(c)}}+O(HAU^{-1})\\
&&\ \ \ \ \ \ +\left(H\min(\sqrt{A}, \frac{1}{|f^{\prime}(a)|})\right) +\left(H\min(\sqrt{A}, \frac{1}{|f^{\prime}(b)|})\right).
\end{eqnarray*}

}

\begin{proof}
This is Lemma 2 on p. 71 of the monograph of A. A. Karatsuba and S. M. Voronin \cite{KV}.
\end{proof}

{\bf Lemma 2.6.} {\it  Suppose $\ell\geq 2$ is  a fixed integer. Then we have the estimate
\begin{eqnarray*}
\sum_{n\leq x} \tau_\ell^2(n)\ll x(\log x)^{\ell^2-1}.
\end{eqnarray*}
}

\begin{proof}
This is a well-known result of the analytic number theory. See, for example,
Theorem 7.6.1 of Ramachandra \cite{Ra}.
\end{proof}

{\bf Lemma 2.7.} {\it Suppose $T\geq 2.$ The asymptotic formula
\begin{equation}\label{Delta-1}
\int_0^T \Delta(x)dx=\frac T4+O(T^{3/4})\end{equation}
holds. Furthermore,  we have the estimate
\begin{eqnarray}\label{mean-delta-2}
&&\int_0^T\Delta^2(x)dx\ll T^{3/2}.
\end{eqnarray}
}

\begin{proof}
The formula (\ref{Delta-1}) is a weaker form of a result of Vorono\"{\i} \cite{Vo}.
The formula (\ref{mean-delta-2}) is a weaker form of the mean square of
$\Delta(x)$ in H. Cram\'er \cite{Cr}.
See, also, J. Furuya \cite{Fu1}.
\end{proof}

{\bf Lemma 2.8.}  {\it   We have the estimate
\begin{eqnarray}
\sum_{n\leq x}|\tau(n)\Delta(n)|\ll x^{5/4}\log^{3/2} x.
\end{eqnarray}}

\begin{proof}
It follows from Lemma 2.6 with $\ell=2$,  Cauchy's inequality and the estimate
$$\sum_{n\leq x}\Delta^2(n)\ll x^{3/2},$$
which is an easy consequence of (2.5).
\end{proof}

{\bf Lemma 2.9.}  {\it   We have the following asymptotic formula
\begin{eqnarray}
\sum_{n\leq x}\tau(n)\Delta(n)= x\sum_{j=0}^3c_j(\log x)^j+O(x^{3/4}\log x),
\end{eqnarray}
where $c_j(j=0, 1, 2, 3)$ are constants.}

\begin{proof}
This is contained in Corollary 1 of Furuya \cite{Fu}.
\end{proof}

{\bf Lemma 2.10.}   {\it  We have the asymptotic formula
$$D_3(x)=C_2x\log^2x+C_1x\log x+C_0x+O(x^{\frac{43}{96}+\varepsilon}),$$
where  $C_2=1/2, C_1, C_0$ are computable constants. }

\begin{proof}
 This is proved in Kolesnik \cite{Ko}.
\end{proof}

{\bf Lemma 2.11.}  {\it Suppose $a_1,\cdots, a_N\in {\Bbb C}$ are arbitrary complex numbers, $T\geq 3$ is a real number. Then we have
$$\int_0^T\left|\sum_{n\leq N}a_nn^{it}\right|^2dt=T\sum_{n\leq N}|a_n|^2
+O\left(\sum_{n\leq N}n|a_n|^2\right).$$
 }

\begin{proof}
 This is Theorem 5.2 of Ivi\'c \cite{Iv}.
\end{proof}

 {\bf Lemma 2.12.}  {\it  Let $F(s)$ be regular in the region $\mathfrak{D}$: $\alpha\leq \sigma\leq \beta, \ t\geq 1,$ and let, for $s\in \mathfrak{D}, F(s)\ll e^{Ct^2}$ for some positive constant $C>0.$ Then for any fixed $q>0$ and $\alpha<\gamma<\beta$ we have
 \begin{eqnarray*}
 \int_1^T|F(\gamma+it)|^qdt\ll \left(\int_1^T|F( \alpha+it)|^qdt+1\right)^{\frac{\beta-\gamma}{\beta-\alpha}}
 \times \left(\int_1^T|F( \beta+it)|^qdt+1\right)^{\frac{ \gamma-\alpha}{\beta-\alpha}}.
 \end{eqnarray*}
 }

\begin{proof}
 This is Lemma 8.3 of Ivi\'c \cite{Iv}.
\end{proof}



{\bf Lemma 2.13.}  {\it Let $A(x)$ be a non-decreasing function such that
the integral
$$H(\xi):=\int_0^\infty e^{-\xi  x}dA(x)$$
converges for all $\xi>0.$ Suppose there exist two real numbers $C\geq 0, \omega>0$
and a non-decreasing function $\psi(x)$ such that
\begin{eqnarray}
\psi(x)\rightarrow \infty, \ \ \psi(x)/x^\omega \mbox{is non-increasing for large}\  x,
\end{eqnarray}
 and
 \begin{eqnarray}
 H(\xi)=\left\{C+O\left(\frac{1}{\psi(1/ \xi)}\right)\right\}
 \xi^{-\omega}\ \  (\xi\rightarrow 0+).
 \end{eqnarray}
Then we have
\begin{eqnarray}
 A(x)=\left\{C+O\left(\frac{1}{\log \psi(x)}\right)\right\}
 \frac{x^{\omega}}{\Gamma(\omega+1)}\ \  (x\rightarrow \infty).
 \end{eqnarray}
}

\begin{proof}
This is Theorem 7.10 in the second part of Tenenbaum \cite{Te}.
\end{proof}

\section{A weighted Piltz divisor problem}

 Define the arithmetic function $\tau_3^{(1)}(n)$   by
  $$\tau_3^{(1)}(n): =\sum_{n=n_1n_2n_3}\log n_3,$$
 which can be viewed as an analogue of $\tau_3(n).$
It is easy to see that
\begin{equation}\label{diri-series}
\sum_{n=1}^\infty\frac{\tau_3^{(1)}(n)}{n^s}=-\zeta^2(s)\zeta^{\prime}(s),\ \ \sigma>1.
\end{equation}
Let
$$D_3^{(1)}(x):=\sum_{n\leq x}\tau_3^{(1)}(n).$$
Then by Perron's formula (see the formula (\ref{perron}) below) it is easy to show that the asymptotic formula
\begin{eqnarray}\label{definition}
D_3^{(1)}(x)=xP_{3}(\log x)+\Delta_3^{(1)}(x)
\end{eqnarray}holds,
where
\begin{eqnarray}\label{residue}
xP_{3}(\log x):&&=Res_{s=1}\left(-\zeta^2(s)\zeta^{\prime}(s)\frac{x^s}{s}\right)\\
&&=x(K_3\log^3 x+K_2\log^2 x+K_1\log  x+K_0)\nonumber
\end{eqnarray}
and the error term $\Delta_3^{(1)}(x)$ satisfies
\begin{equation}\label{bound-1}
 \Delta_3^{(1)}(x)=O(x^{1/2+\varepsilon}).
\end{equation} Here $K_j(j=0, 1, 2,3)$ are computable constants.
Note that $x^\varepsilon$ in (3.4) can be replaced by $\log^4 x.$

 In this section we shall prove a Vorono\"{\i} formula of
 $ \Delta_3^{(1)}(x)$ and  give an improvement of (\ref{bound-1}).

\subsection{A Vorono\"{\i} formula of $ \Delta_3^{(1)}(x)$}\

Let $1\ll N\ll x$ be an integer. Suppose $1\ll T\ll x$ is a parameter to be determined.  Actually, $T$ will be determined by the formula (\ref{T}) below.  By Perron's formula and (\ref{diri-series}) we get
\begin{equation}\label{perron}
D_3^{(1)}(x)=\frac{1}{2\pi i}\int_{1+\varepsilon-iT}^{1+\varepsilon+iT}-\zeta^2(s)\zeta^{\prime}(s)\frac{x^s}{s}ds
+O\left(\frac{x^{1+\varepsilon}}{T}\right).
\end{equation}

By the residue theorem we have
\begin{equation}\label{shift}
\Delta_3^{(1)}(x)=I_1(x)+I_2(x)-I_3(x)+O\left(\frac{x^{1+\varepsilon}}{T}\right),
\end{equation}
where
\begin{eqnarray*}
&&I_1(x):=\frac{1}{2\pi i}\int_{-\varepsilon+iT}^{1+\varepsilon+iT}-\zeta^2(s)\zeta^{\prime}(s)\frac{x^s}{s}ds,\\
&&I_2(x):=\frac{1}{2\pi i}\int_{-\varepsilon-iT}^{-\varepsilon+iT}-\zeta^2(s)\zeta^{\prime}(s)\frac{x^s}{s}ds,\\
&&I_3(x):=\frac{1}{2\pi i}\int_{-\varepsilon-iT}^{1+\varepsilon-iT}-\zeta^2(s)\zeta^{\prime}(s)\frac{x^s}{s}ds.
\end{eqnarray*}
Note that the residue $\zeta^2(0)\zeta^{\prime}(0)$ at $s=0$ was absorbed into the error term
$O(x^{1+\varepsilon}/T).$

By Lemma 2.3 it is easy to see that
\begin{eqnarray}\label{horizontal}
 I_1(x) \ll T^{-1}(x+T^{3/2})\log^4 T, \ \ \ \ \ I_3(x) \ll T^{-1}(x+T^{3/2})\log^4 T.
\end{eqnarray}

It suffices for us to evaluate $I_2(x).$ By Lemma 2.1 we have
\begin{eqnarray}\label{left}
I_2(x)&&
=\frac{1}{2\pi i}\int_{-\varepsilon-iT}^{-\varepsilon+iT}\zeta^2(s)
\left(\chi(s)\zeta^{\prime}(1-s)-\chi^{\prime}(s)\zeta(1-s) \right)\frac{x^s}{s}ds\nonumber\\
&&=\frac{1}{2\pi i}\int_{-\varepsilon-iT}^{-\varepsilon+iT}
 (-1)\chi^2(s)\chi^{\prime}(s)\zeta^3(1-s)  \frac{x^s}{s}ds\nonumber\\
&&\ \ +
\frac{1}{2\pi i}\int_{-\varepsilon-iT}^{-\varepsilon+iT}
 \chi^3(s)\zeta^2(1-s)\zeta^{\prime}(1-s) \frac{x^s}{s}ds\nonumber\\
 &&=I_{21}(x)+I_{22}(x),\nonumber
\end{eqnarray}
say. We will evaluate  $I_{21}(x)$ only. The term  $I_{22}(x)$ can be evaluated similarly.

  It is easy to see that
\begin{eqnarray}\label{I21}
I_{21}(x)&& = \frac{1}{2\pi i}\int_{-\varepsilon-iT}^{-\varepsilon+iT}
 (-1)\chi^2(s)\chi^{\prime}(s)\zeta^3(1-s)  \frac{x^s}{s}ds
 =2\Re I+O(x^{-\varepsilon})
 \end{eqnarray}
with
\begin{equation}
I:=\frac{1}{2\pi}\label{I}
\int_1^T (-1)\chi^2(-\varepsilon+it)\chi^{\prime}(-\varepsilon+it)
\sum_{n=1}^\infty \frac{\tau_3(n)n^{it}}{n^{1+\varepsilon}}  \frac{x^{-\varepsilon+it}}{-\varepsilon+it}dt,
\end{equation}
where the infinite series in (\ref{I}) is absolutely convergent. So after changing  the order of integration and summation we get
\begin{eqnarray}\label{I-form}
I=\frac{x^{-\varepsilon}}{2\pi}\sum_{n=1}^\infty \frac{\tau_3(n) }{n^{1+\varepsilon}}
\int_1^T (-1)\chi^2(-\varepsilon+it)\chi^{\prime}(-\varepsilon+it)
   \frac{(nx)^{it}}{-\varepsilon+it}dt.
\end{eqnarray}

Using Lemma 2.2 we have
\begin{eqnarray}\label{function}
&& \ \ \ \ \ -\chi^2(-\varepsilon+it)\chi^{\prime}(-\varepsilon+it)
   \frac{(nx)^{it}}{-\varepsilon+it}\\
&&=-
\chi^3(-\varepsilon+it)\frac{\chi^{\prime}(- \varepsilon+it)}{\chi(-\varepsilon+it)}
\frac{(nx)^{ it}}{-\varepsilon+it}\nonumber\\
&&=-\left(\frac{t}{2\pi}\right)^{ \frac 32+3\varepsilon-3it}e^{3i(t+\frac \pi 4)}
\frac{(nx)^{it}}{it}
\left(1+O\left(\frac 1t\right)\right)
\left(-\log\frac{t}{2\pi}+O\left(\frac 1t\right)\right)\nonumber\\
&&= \frac{e^{\frac{\pi i}{4}}}{2\pi}
 \left(\log\frac{t}{2\pi}\right)
\left(\frac{t}{2\pi}\right)^{\frac 12+ 3 \varepsilon}  e(G(t))
  +O(t^{-\frac 12 +3 \varepsilon} \log t),\nonumber
\end{eqnarray}
where $$2\pi G(t)=2\pi G(t;n,x):=-3t\log\frac{t}{2\pi}+3t+t\log nx.$$

From (\ref{I-form}) and (\ref{function}) we get
\begin{eqnarray}\label{I-newform}
I&&=\frac{x^{-\varepsilon}e^{\frac{i\pi}{4}}}{(2\pi)^2}\sum_{n=1}^\infty \frac{\tau_3(n) }{n^{1+\varepsilon}}
I(n,x,T)  +O(x^{-\varepsilon}T^{\frac 12+3\varepsilon}),
\end{eqnarray}
where
$$I(n,x,T):=\int_1^T \left(\log\frac{t}{2\pi}\right)
\left(\frac{t}{2\pi}\right)^{\frac 12+ 3 \varepsilon}  e(G(t)) dt.$$

We choose the parameter $T$ satisfying
\begin{equation}\label{T}
 N+\frac 12=\left(\frac{T}{2\pi }\right)^3\frac 1x, \ \ \ N\in {\Bbb N}.
\end{equation}
The condition $1\ll N\ll x $ ensures that
$x^{1/3}\ll T\ll x^{2/3}.$

Let
\begin{eqnarray}\label{I1}
I_{1}:&&=\frac{x^{-\varepsilon}e^{\frac{i\pi}{4}}}{(2\pi)^2}\sum_{n\leq N} \frac{\tau_3(n) }{n^{1+\varepsilon}}
I(n,x,T),\\
I_{2}:&&=\frac{x^{-\varepsilon}e^{\frac{i\pi}{4}}}{(2\pi)^2}\sum_{n>N} \frac{\tau_3(n) }{n^{1+\varepsilon}}
 I(n,x,T).\nonumber
\end{eqnarray}

It is easy to see that $G^{\prime}(t)=\frac{1}{2\pi}\log\frac{(2\pi)^3nx}{t^3}.$
If $1\leq t\leq T, n>N,$ then we have
$$G^{\prime}(t)=\frac{1}{2\pi}\log\frac{nx}{(t/2\pi)^3}\geq \frac{1}{2\pi}\log \frac{n}{N+1/2}.$$

So by partial integration and the first derivative test we have
\begin{eqnarray}\label{I2-upperbound}
I_{2}&&\ll x^{-\varepsilon}T^{\frac 12+4\varepsilon}
\left(\sum_{N<n\leq 3N/2} \frac{\tau_3(n) }{n^{1+\varepsilon}}
\frac{1}{ \log \frac{n}{N+1/2}}
+\sum_{n>3N/2} \frac{\tau_3(n) }{n^{1+\varepsilon}}\right)\\
&&\ll T^{\frac 12}x^{\varepsilon}, \nonumber
\end{eqnarray}
where we write $n=N+r( 1\leq r\leq N/2)$ in the first sum in (\ref{I2-upperbound}).

Now we consider $I(n,x,T)$ in $I_1$.
The saddle point is $t_0:=2\pi(nx)^{1/3},$ namely, $G^{\prime}(t_0)=0.$
We have $t_0\leq T$ for any $n\leq N.$
It is easy to see that
\begin{eqnarray*}
&&G(t_0)=3(nx)^{1/3},\ \ G^{\prime\prime}(t_0)=-\frac{3}{(2\pi)^2 (nx)^{1/3}},\\
&&\left(\log\frac{t_0}{2\pi}\right)
\left(\frac{t_0}{2\pi}\right)^{\frac 12+ 3 \varepsilon}
=\frac{1}{3} (nx)^{\frac 16+\varepsilon}\log (nx).
\end{eqnarray*}

We will use Lemma 2.5 first with
$$g(x)=\left(\log\frac{x}{2\pi}\right)\left(\frac{x}{2\pi}\right)^{\frac 12+ 3 \varepsilon},\ \ f(x)=G(x),$$
and then use a conjugate procedure.
We shall show that
 \begin{eqnarray}\label{I(n,x,T)-formula}
I(n,x,T)
&&=\frac{2\pi e^{-\frac{\pi i}{4}}}{3\sqrt 3}(nx)^{\frac 13+\varepsilon}\log (nx)
\times e(3(nx)^{1/3})\\
&&\ \ \ \
  + O( T^{\frac 12+3\varepsilon}\log T)+O\left( T^{\frac 12+3\varepsilon}\log T\times \frac{1}{\log\frac{n}{N+1/2}}\right).   \nonumber
\end{eqnarray}

{\bf Case I: } $2t_0<T.$
In this case we write
\begin{eqnarray}\label{sum}
&&  I(n,x,T)=
\int_1+\int_2+\int_3,
\end{eqnarray}
where
\begin{eqnarray}\label{int1}
&&\int_1:=\int_1^{t_0/2} \left(\log\frac{t}{2\pi}\right)
\left(\frac{t}{2\pi}\right)^{\frac 12+ 3 \varepsilon}  e(G(t)) dt,\\
&&\int_2:=\int_{t_0/2}^{2t_0} \left(\log\frac{t}{2\pi}\right)
\left(\frac{t}{2\pi}\right)^{\frac 12+ 3 \varepsilon}  e(G(t)) dt,\nonumber\\
&&\int_3:=\int_{2t_0}^T \left(\log\frac{t}{2\pi}\right)
\left(\frac{t}{2\pi}\right)^{\frac 12+ 3 \varepsilon}  e(G(t)) dt.\nonumber
 \end{eqnarray}

By partial summation and the first derivative test we
have
 \begin{equation}\label{upperbound-int1}
 \int_1\ll T^{\frac 12+3\varepsilon}\log T,\ \ \int_3\ll T^{\frac 12+3\varepsilon}\log T.
 \end{equation}

From Lemma 2.5 we get that
\begin{eqnarray}\label{use-Lemma}
&&\ \ \ \ \ \ \ \int_{t_0/2}^{2t_0} \left(\log\frac{t}{2\pi}\right)
\left(\frac{t}{2\pi}\right)^{\frac 12+ 3 \varepsilon}  e(G(t)) dt\\
&&
=\frac{2\pi e^{-\frac{\pi i}{4}}}{3\sqrt 3}(nx)^{\frac 13+\varepsilon}\log (nx)
\times e(3(nx)^{1/3})
  + O(T^{\frac 12+3\varepsilon}\log T),  \nonumber
\end{eqnarray}
where we note that $|G^{\prime}(t_0/2)|\gg 1, |G^{\prime}(2t_0)|\gg 1.$

From (\ref{sum})-(\ref{use-Lemma}) we get (\ref{I(n,x,T)-formula}) for  $2t_0<T.$

{\bf Case II:}  $T\leq 2t_0.$
In this case we write
\begin{eqnarray}\label{sum2}
&&  I(n,x,T)=
\int_1+ \int_4,
\end{eqnarray}
where
\begin{eqnarray*}
&&\int_4:=\int_{t_0/2}^T \left(\log\frac{t}{2\pi}\right)
\left(\frac{t}{2\pi}\right)^{\frac 12+ 3 \varepsilon}  e(G(t)) dt
 \end{eqnarray*}
and where $\int_1$ was defined in (\ref{int1}).

By Lemma 2.5 again we have
\begin{eqnarray}\label{use2-Lemma}
&&\ \ \ \ \ \ \ \int_{t_0/2}^T \left(\log\frac{t}{2\pi}\right)
\left(\frac{t}{2\pi}\right)^{\frac 12+ 3 \varepsilon}  e(G(t)) dt\\
&&
=\frac{2\pi e^{\frac{-\pi i}{4}}}{3\sqrt 3}(nx)^{\frac 13+\varepsilon}\log (nx)
\times e(3(nx)^{1/3})\nonumber\\&&\ \ \
  + O(T^{\frac 12+3\varepsilon}\log T )+O\left(T^{\frac 12+3\varepsilon}\log T\times \frac{1}{|\log\frac{n}{N+1/2}|}\right).   \nonumber
\end{eqnarray}

Now (\ref{I(n,x,T)-formula}) follows from (\ref{upperbound-int1}),
 (\ref{sum2})  and (\ref{use2-Lemma}).

The contribution of the two error terms in (\ref{I(n,x,T)-formula} ) to $I_1$ is
\begin{eqnarray*}
&& \ll  x^{-\varepsilon}T^{\frac 12+3\varepsilon}\log T \sum_{n\leq N} \frac{\tau_3(n) }{n^{1+\varepsilon}}+x^{-\varepsilon}T^{\frac 12+3\varepsilon}\log T \sum_{n\leq N} \frac{\tau_3(n) }{n^{1+\varepsilon}} \frac{1}{|\log\frac{n}{N+1/2}|}\\
&& \ll x^{-\varepsilon}T^{\frac 12+3\varepsilon}\log T \sum_{n\leq N/2} \frac{\tau_3(n) }{n^{1+\varepsilon}} \frac{1}{|\log\frac{n}{N+1/2}|}\nonumber\\&&\ \ \
+x^{-\varepsilon}T^{\frac 12+3\varepsilon}\log T \sum_{N/2<n\leq N} \frac{\tau_3(n) }{n^{1+\varepsilon}} \frac{1}{|\log\frac{n}{N+1/2}|}\nonumber\\
&&\ll x^{-\varepsilon}T^{\frac 12+3\varepsilon}\log T,\nonumber
\end{eqnarray*}
where we used $|\log \frac{n}{N+1/2}|\gg 1$ for $1\leq n\leq N/2$ and we wrote
$n=N-r$ for $N/2<n< N.$ The contribution of the main term in (\ref{I(n,x,T)-formula} ) to $I_1$ is
\begin{eqnarray*}
&&=\frac{x^{-\varepsilon}e^{\frac{i\pi}{4}}}{(2\pi)^2}\sum_{n\leq N} \frac{\tau_3(n) }{n^{1+\varepsilon}}\frac{2\pi e^{\frac{\pi i}{4}}}{3\sqrt 3}(nx)^{\frac 13+\varepsilon}\log (nx)
\times e(3(nx)^{1/3})\\
&&=\frac{x^{\frac 13}}{6\sqrt{3}\pi}\sum_{n\leq N} \frac{\tau_3(n)\log(nx) }{n^{\frac 23}}e\left(3(nx)^{1/3} \right)\nonumber
\end{eqnarray*}

Thus from the above two estimates we get
\begin{eqnarray}\label{I1-finally}
&&I_1=\frac{x^{\frac 13}}{6\sqrt{3}\pi}\sum_{n\leq N} \frac{\tau_3(n)\log(nx) }{n^{\frac 23}}e\left(3(nx)^{1/3} \right) +O(x^{-\varepsilon}T^{\frac 12+3\varepsilon}\log T).
\end{eqnarray}
From (\ref{I-newform}), (\ref{I1}), (\ref{I2-upperbound}) and (\ref{I1-finally}) we get
\begin{eqnarray}\label{I21}
I_{21}(x)=\frac{x^{\frac 13}}{3\sqrt{3}\pi}\sum_{n\leq N}
\frac{\tau_3(n)\log(nx) }{n^{\frac 23}}
  \cos(6\pi (nx)^{1/3})
+O(x^{\frac 23+\varepsilon}N^{-\frac 13 } ).
\end{eqnarray}

Similarly we can get
\begin{eqnarray}\label{I22}
I_{22}(x)=-\frac{x^{1/3}}{\sqrt{3}\pi}\sum_{n\leq N}\frac{\tau_3^{(1)}(n)}{n^{2/3}}\cos(6\pi (nx)^{1/3})
+O(x^{\frac 23+\varepsilon}N^{-\frac 13 } ).
\end{eqnarray}

From (\ref{shift})-(\ref{left}), (\ref{I21}) and  (\ref{I22}), we get the following proposition.

{\bf Proposition 3.1.} {\it Suppose $N\in {\Bbb N}$ with  $1\ll N\ll x.$ Then we have
\begin{eqnarray}\label{VV1}
\Delta_3^{(1)}(x)&&=\frac{x^{1/3}}{3\sqrt{3}\pi}\sum_{n\leq N}\frac{\tau_3(n)\log (nx)}{n^{2/3}}\cos(6\pi (nx)^{1/3})\\
&&\ \  -\frac{x^{1/3}}{\sqrt{3}\pi}\sum_{n\leq N}\frac{\tau_3^{(1)}(n)}{n^{2/3}}\cos(6\pi (nx)^{1/3})
+O(x^{2/3+\varepsilon}N^{-1/3}).\nonumber
\end{eqnarray} }

\subsection{An upper bound of $\Delta_3^{(1)}(x)$}\

For $\Delta_3(x),$ we have the similar Voronoi's formula. Suppose $N\in {\Bbb N}$ with  $1\ll N\ll x.$ Then we have
\begin{eqnarray}\label{VVV2}
\Delta_3(x)=  \frac{x^{1/3}}{\sqrt{3}\pi}\sum_{n\leq N}\frac{\tau_3(n)}{n^{2/3}}\cos(6\pi (nx)^{1/3})
  +O(x^{2/3+\varepsilon}N^{-1/3}).
\end{eqnarray} The above formula is a special case $k=3$ of Lemma 3 in
  Ivi\'{c} and Zhai \cite{IZ}.

Comparing the Vorono\"{\i} formula (\ref{VV1}) of $\Delta_3^{(1)}(x)$ and
the Vorono\"{\i} formula (\ref{VVV2}) of $\Delta_3(x)$, we find that they are of  the same structure. Choosing $N=x^{21/32}$ in (\ref{VVV2}), then the error term
$O(x^{2/3+\varepsilon}N^{-1/3})$  is now $O(x^{43/96+\varepsilon}).$
Kolesnik \cite{Ko} proved that
\begin{eqnarray}\label{ex1}
  \frac{x^{1/3}}{\sqrt{3}\pi}\sum_{n\leq N}\frac{\tau_3(n)}{n^{2/3}}\cos(6\pi (nx)^{1/3})
 \ll x^{43/96+\varepsilon},
\end{eqnarray}
or more precisely, he proved that the estimate
\begin{eqnarray}\label{ex22}
  \sum_{n_1\sim N_1, n_2\sim N_2, n_3\sim N_3 }
  \frac{ 1}{(n_1n_2n_3)^{2/3}}\cos(6\pi (n_1n_2n_3x)^{1/3})
 \ll x^{11/96+\varepsilon}
\end{eqnarray}
holds, where $N_1N_2N_3\ll x^{21/32}.$

Choosing $N=x^{21/32}$ in (\ref{VV1}), then the error term
$O(x^{2/3+\varepsilon}N^{-1/3})$  is   $O(x^{43/96+\varepsilon})$ again.
By partial integration and then using (\ref{ex1}), we find that the first sum
in (\ref{VV1})  is bounded by
$ x^{43/96+\varepsilon}.$ Using a splitting argument to the second sum in
(\ref{VV1}), we only need to estimate sums of the form
\begin{eqnarray*}
  \sum_{n_1\sim N_1, n_2\sim N_2, n_3\sim N_3 }
  \frac{\log n_3}{(n_1n_2n_3)^{2/3}}\cos(6\pi (n_1n_2n_3x)^{1/3}).
\end{eqnarray*}
Using a partial summation to remove the factor $\log n_3$ in the above sum, and then using
(\ref{ex22}) we get that the second sum in (\ref{VV1}) is also bounded by
$ x^{43/96+\varepsilon}.$ Collecting all estimates we get the following

{\bf Proposition 3.2.}  {\it The estimate
\begin{eqnarray*}
\Delta_3^{(1)}(x)\ll x^{43/96+\varepsilon}
\end{eqnarray*}holds.
}

\section{A variant of the Dirichlet divisor problem}

Consider the sum
\begin{equation}\label{4.1}
U(x):=\sum_{\stackrel{m,n\leq x}{gcd(m,n)=1}}\tau(m)\tau(n),
\end{equation}
which can be viewed as a generalization of $D(x).$ The sum
$U(x)$ plays an important role when studying $S^{(2)}(x)$ and $C^{(2)}(x)$.
So it has its own interest.

Nowak and Toth \cite{NT} proved that $U(x)$ has the asymptotic formula
\begin{equation}\label{4.2}
U(x)=x^2(b_2\log^2 x+b_1\log x+b_0)+E_U(x),
\end{equation}
where $b_0, b_1, b_2$ are real constants defined by (19) therein and the error term
$E_U(x)$ in (\ref{4.2}) satisfies
\begin{equation}\label{4.3}
E_U(x)= O(x^{4/3+\varepsilon}).
\end{equation}

In this section
we prove more results about $E_U(x).$  We will prove the following three propositions.

{\bf Proposition 4.1.}  {\it  Suppose $1/4<\theta<1/3$ is a real number such that
$\Delta(t)\ll t^\theta$. Then  the estimate
\begin{equation}\label{4.4}
E_U(x)\ll x^{1+\theta}\log x
\end{equation}
holds.}

{\bf Proposition 4.2.}  {\it
  Let $T\geq 3$ be a real parameter, then
\begin{equation}\label{4.5}
\int_1^T E_U(x)dx\ll T^{2+\varepsilon}.
\end{equation}}

{\bf Proposition 4.3.}  {\it
   Let $T\geq 3$ be a real parameter, then we have  the estimate
\begin{equation}\label{4.6}
\int_1^T(E_U(x))^2dx\ll
 T^{\frac 72 }\log^2 T.
\end{equation}
}

\subsection{An expression of $E_U(x)$}\

In this subsection we give an expression of $E_U(x)$.

For $s=\sigma+it\in {\Bbb C}$ and $ z=u+iv\in {\Bbb C}$ with $ \sigma>1$ and $ u>1$, we define the Dirichlet series
\begin{equation*}
\mathfrak{D}(s,z):=\sum_{\stackrel{m=1,n=1}{gcd(m,n)=1}}^\infty\frac{\tau(m)\tau(n)}{m^sn^z}.
\end{equation*}

It is easy to see that  we have the Euler product form
\begin{eqnarray}\label{4.7}
\mathfrak{D}(s,z)&&= \prod_{p}\left(\sum_{\stackrel{\alpha,\beta\geq 0}{\alpha\beta=0}}^\infty
\frac{(\alpha+1)(\beta+1)}{p^{\alpha s+\beta z}}\right)\\
&&=\prod_{p}\left(1+\sum_{\alpha=1}^\infty
\frac{\alpha+1 }{p^{\alpha s }}+\sum_{\beta=1}^\infty
\frac{\beta+1 }{p^{\beta z }}
\right).\nonumber
\end{eqnarray}

By the expression $(|\eta|<1)$
$$(1-\eta)^{-2}-1=\sum_{n=1}^\infty (n+1)\eta^n$$ we have
 \begin{eqnarray}\label{4.8}
&&\ \ \ \  1+\sum_{\alpha=1}^\infty
\frac{\alpha+1 }{p^{\alpha s }}+\sum_{\beta=1}^\infty
\frac{\beta+1 }{p^{\beta z }}=
\frac{1}{(1-p^{-s})^2}+\frac{1}{(1-p^{-z})^2}-1 \\
&&=\frac{(1-p^{-s})^2+(1-p^{-z})^2 -(1-p^{-s})^2(1-p^{-z})^2  }{(1-p^{-s})^2(1-p^{-z})^2}\nonumber\\
&&=\frac{1-4p^{-s-z}
+2p^{-2s-z}+2p^{-s-2z}-p^{-2s-2z} }{(1-p^{-s})^2(1-p^{-z})^2}.\nonumber
\end{eqnarray}

From (\ref{4.7}) and (\ref{4.8}) we get the following Lemma 4.1.

{\bf Lemma 4.1.} {\it The Dirichlet series
   $\mathfrak{D}(s,z)$ can be written as the form
\begin{equation}\label{4.9}
\mathfrak{D}(s,z)=\zeta^2(s)\zeta^2(z) \mathfrak{D}_1(s,z),
\end{equation}
where $\mathfrak{D}_1(s,z)$ can be written  as a Dirichlet series
$$  \mathfrak{D}_1(s,z)=\sum_{m=1}^\infty\sum_{n=1}^\infty\frac{d_1(m,n)}{m^sn^z},
$$
 which is absolutely convergent for $\sigma+u>1.$ Specifically for any fixed integers $\ell_1\geq 1$ and $\ell_2\geq 1,$ the estimate
 \begin{equation}\label{4.10}
  \mathfrak{D}_1(s,z)=\sum_{m=1}^
  \infty\sum_{n=1}^\infty\frac{|d_1(m,n)|\tau_{\ell_1}(m)\tau_{\ell_2}(n)}
  {m^\sigma n^u}\ll_\varepsilon 1
  \end{equation}
holds uniformly for $ \sigma\geq \varepsilon, \ \ u\geq \varepsilon,\ \ \sigma+u\geq 1+\varepsilon.$
 }

From (\ref{4.9})   we see that if $gcd(m,n)=1,$ then
$$\tau(m)\tau(n)=\sum_{\stackrel{m=m_1m_2}{n=n_1n_2}}d_1(m_1,n_1)\tau(m_2)\tau(n_2),$$
which, with the help of (\ref{1.1}), gives
\begin{eqnarray}\label{4.11}
U(x)&&=\sum_{\stackrel{m_1\leq x}{n_1\leq x}}d_1(m_1,n_1)
\sum_{m_2\leq \frac{x}{m_1}}\tau(m_2)
\sum_{n_2\leq\frac{x}{n_1}}\tau(n_2)\\
&&=\sum_{\stackrel{m_1\leq x}{n_1\leq x}}d_1(m_1,n_1)
\left(M\left(\frac{x}{m_1}\right)+\Delta\left(\frac{x}{m_1}\right)\right)
\nonumber\\ &&\ \ \ \ \ \ \ \ \ \ \ \times
\left(M\left(\frac{x}{n_1}\right)+\Delta\left(\frac{x}{n_1}\right)\right)\nonumber\\
&&=U_1(x)+U_2(x)+U_3(x)+U_4(x),\nonumber
\end{eqnarray}say,
where
\begin{eqnarray*}
U_1(x)&&=\sum_{\stackrel{m_1\leq x}{n_1\leq x}}d_1(m_1,n_1)
 M\left(\frac{x}{m_1}\right)M\left(\frac{x}{n_1}\right),\\
 U_2(x)&&=\sum_{\stackrel{m_1\leq x}{n_1\leq x}}d_1(m_1,n_1)
 M\left(\frac{x}{m_1}\right)\Delta\left(\frac{x}{n_1}\right),\\
 U_3(x)&&=\sum_{\stackrel{m_1\leq x}{n_1\leq x}}d_1(m_1,n_1)
  \Delta\left(\frac{x}{m_1}\right)M\left(\frac{x}{n_1}\right),\\
 U_4(x)&&=\sum_{\stackrel{m_1\leq x}{n_1\leq x}}d_1(m_1,n_1)
  \Delta\left(\frac{x}{m_1}\right)\Delta\left(\frac{x}{n_1}\right),\\
  M(x)&&=x\log x+(2\gamma-1)x.
\end{eqnarray*}

By some easy calculations with the help of (\ref{4.10}) we get
\begin{eqnarray}\label{4.12}
U_1(x)=x^2(b_2\log^2 x+b_1\log x+b_0)+O(x^{1+\varepsilon}),
\end{eqnarray}
where $b_j(j=0,1,2)$ are defined in (\ref{4.2}).
Using the estimate $\Delta(t)\ll t^{1/2+\varepsilon/2}$ we get
\begin{eqnarray} \label{4.13}
 U_4(x)&&\ll\sum_{\stackrel{m_1\leq x}{n_1\leq x}}|d_1(m_1,n_1)|
   \left(\frac{x}{m_1}\right)^{1/2+\varepsilon/2} \left(\frac{x}{n_1}\right)^{1/2+\varepsilon/2}\\
   &&\ll x^{1+\varepsilon}\sum_{\stackrel{m_1\leq x}{n_1\leq x}}
   \frac{|d_1(m_1,n_1)|}{(m_1n_1)^{1/2+\varepsilon/2}}\ll  x^{1+\varepsilon},\nonumber
\end{eqnarray}
where in the last step we used   (\ref{4.10}) of Lemma 4.1.

By symmetry we have $U_2(x)=U_3(x),$ which combining with (\ref{4.11})-(\ref{4.13}) gives

{\bf Lemma 4.2.}  {\it
We have  the expression
$$E_U(x)=2\sum_{\stackrel{m_1\leq x}{n_1\leq x}}d_1(m_1,n_1)
 M\left(\frac{x}{m_1}\right)\Delta\left(\frac{x}{n_1}\right)+O(x^{1+\varepsilon}).$$
}

\subsection{Proof of Proposition 4.1}\

By using the estimate $\Delta(x)\ll x^{\theta}$ and (\ref{4.10}) we have
\begin{eqnarray*}
E_U(x)&&\ll\sum_{\stackrel{m_1\leq x}{n_1\leq x}}|d_1(m_1,n_1)|
  \frac{x\log x}{m_1}  \left(\frac{x}{n_1}\right)^{\theta}+x^{1+\varepsilon}  \\
  &&\ll x^{1+\theta}\log x
   \times\sum_{\stackrel{m_1\leq x}{n_1\leq x}}\frac{|d_1(m_1,n_1)|}{m_1n_1^\theta} +x^{1+\varepsilon}\nonumber\\
   &&\ll x^{1+\theta}\log x.\nonumber
\end{eqnarray*}

\subsection{Proof of Proposition 4.2}\

Let
$$E_U^{*}(x)=2\sum_{\stackrel{m_1\leq x}{n_1\leq x}}d_1(m_1,n_1)
 M\left(\frac{x}{m_1}\right)\Delta\left(\frac{x}{n_1}\right).$$

By using the expression
$$M\left(\frac{x}{m_1}\right)=\frac{x}{m_1}(\log x-\log m_1+2\gamma-1)$$
we can write
\begin{eqnarray}\label{4.14}
&&E_U^{*}(x)=E_{U1}^{*}(x)+E_{U2}^{*}(x)+E_{U3}^{*}(x),\\
&&E_{U1}^{*}(x)=
2\sum_{\stackrel{m_1\leq x}{n_1\leq x}}d_1(m_1,n_1)
  \frac{x\log x}{m_1} \Delta\left(\frac{x}{n_1}\right),\nonumber\\
&&E_{U2}^{*}(x)=-
2\sum_{\stackrel{m_1\leq x}{n_1\leq x}}d_1(m_1,n_1)
  \frac{x\log m_1}{m_1} \Delta\left(\frac{x}{n_1}\right),\nonumber\\
  &&E_{U3}^{*}(x)=
2(2\gamma-1)\sum_{\stackrel{m_1\leq x}{n_1\leq x}}d_1(m_1,n_1)
  \frac{x }{m_1} \Delta\left(\frac{x}{n_1}\right).\nonumber
\end{eqnarray}

It is easy to see that
\begin{eqnarray}\label{4.15}
\int_1^TE_{U1}^{*}(x)dx=2\sum_{\stackrel{m_1\leq T}{n_1\leq T}}
  \frac{ d_1(m_1,n_1)}{m_1}\int_{\max(m_1,n_1)}^T x\log x  \Delta\left(\frac{x}{n_1}\right)dx.
\end{eqnarray}
We have
\begin{eqnarray}\label{4.16}
&&\ \ \int_{\max(m_1,n_1)}^T x\log x  \Delta\left(\frac{x}{n_1}\right)dx\\
&&=n_1^2\int_{\max(m_1,n_1)}^T \frac{x}{n_1}\log x  \Delta\left(\frac{x}{n_1}\right)d\frac{x}{n_1}\nonumber\\
&&=n_1^2\int_{\max(m_1,n_1)}^T \frac{x}{n_1}\left(\log \frac{x}{n_1}+\log n_1 \right) \Delta\left(\frac{x}{n_1}\right)d\frac{x}{n_1}\nonumber\\
&&=n_1^2\int_{\max(m_1/n_1,1)}^{T/n_1} w\left(\log w+\log n_1 \right) \Delta\left(w\right)dw.\nonumber
\end{eqnarray}
We write
$$F(t)=\int_0^t\Delta(w)dw=\frac t4+P(t).$$
By (2.4) of Lemma 2.7 we have $P(t)\ll t^{3/4}.$ Hence
$$ \Delta(w)dw=\frac 14 dw+dP(w).$$
Thus we can write
\begin{eqnarray}\label{4.17}
  n_1^2\int_{\max(m_1/n_1,1)}^{T/n_1} w\left(\log w+\log n_1 \right) \Delta\left(w\right)dw=W_1+W_2,
\end{eqnarray}
where
\begin{eqnarray*}
 W_1&&=\frac{n_1^2}{4}\int_{\max(m_1/n_1,1)}^{T/n_1} w\left(\log w+\log n_1 \right)dw,\\
 W_2&&=n_1^2\int_{\max(m_1/n_1,1)}^{T/n_1} w\left(\log w+\log n_1 \right) dP(w).
\end{eqnarray*}
By partial integration it is easy to check that
\begin{eqnarray}\label{4.18}
 W_1\ll T^2\log T,\ \ \
 W_2 \ll n_1^2 (T/n_1)^{7/4}\log T\ll T^2\log T\ \ (n_1\ll T).
\end{eqnarray}

From (\ref{4.15})-(\ref{4.18}) we get
\begin{eqnarray}\label{4.19}
\int_1^TE_{U1}^{*}(x)dx\ll T^{2}\log T.
\end{eqnarray}

Similarly, we have
\begin{eqnarray}\label{4.20}
&&\int_1^TE_{U2}^{*}(x)dx\ll T^{2}\log T,\ \ \
\int_1^TE_{U3}^{*}(x)dx\ll T^{2}.
\end{eqnarray}

From (\ref{4.10}), (\ref{4.14}), (\ref{4.19}) and (\ref{4.20}) we get
\begin{eqnarray*}
\int_1^TE_U^{*}(x)dx&&\ll T^2\log T\times \sum_{\stackrel{m_1\leq T}{n_1\leq T}}\frac{|d_1(m_1,n_1)|}{m_1} \\
&&\ll T^2\log T\times \sum_{\stackrel{m_1\leq T}{n_1\leq T}}\frac{|d_1(m_1,n_1)|n_1^\varepsilon}{m_1n_1^\varepsilon}\nonumber\\
&&\ll T^{2+\varepsilon},\nonumber
\end{eqnarray*}
which combining with Lemma 4.2 implies Proposition 4.2.

\subsection{Proof of Proposition 4.3}\

 Lemma 4.2, we only need to prove
 that
 \begin{eqnarray}
  \int_1^{T}(E_U^{*}(x))^2dx\ll T^{7/2}\log^2  T.
\end{eqnarray}

Suppose $X\geq 3$ is a  large parameter. We shall bound the integral
$\int_{X/2}^{X}(E_U^{*}(x))^2dx.$
For any $X/2\leq x\leq X,$ by the definition of $E_U^{*}(x)$ we have
 $$E_U^{*}(x)\ll X\log X \sum_{\stackrel{m_1\leq X}{n_1\leq X}}
\frac{| d_1(m_1,n_1)
|}{m_1}\left|\Delta\left(\frac{x}{n_1}\right)\right|.$$

By Cauchy's inequality and (4.10) we have
\begin{eqnarray*}
 (E_U^{*}(x))^2&&\ll X^2\log^2 X\sum_{\stackrel{m_1\leq X}{n_1\leq X}}
\frac{ |d_1(m_1,n_1)|}{m_1n_1^{1/4}}
\sum_{\stackrel{m_1\leq X}{n_1\leq X}}
\frac{ |d_1(m_1,n_1)|n_1^{1/4}}{m_1}\left|\Delta\left(\frac{x}{n_1}\right)\right|^2\\
&&\ll X^2\log^2  X
\sum_{\stackrel{m_1\leq X}{n_1\leq X}}
\frac{| d_1(m_1,n_1)|n_1^{1/4}}{m_1}\left|\Delta\left(\frac{x}{n_1}\right)\right|^2.
\end{eqnarray*}

So by the second assertion of Lemma 2.7 and (4.10) again we have
\begin{eqnarray*}
\int_{X/2}^{X} (E_U^{*}(x))^2dx&&\ll   X^2\log^2  X
\sum_{\stackrel{m_1\leq X}{n_1\leq X}}
\frac{ |d_1(m_1,n_1)|n_1^{1/4}}{m_1}\int_{X/2}^{X}
\left|\Delta\left(\frac{x}{n_1}\right)\right|^2dx\\
&&\ll  X^2\log^2  X
\sum_{\stackrel{m_1\leq X}{n_1\leq X}}
\frac{ |d_1(m_1,n_1)|n_1^{5/4}}{m_1}\int_{X/2}^{X}
\left|\Delta\left(\frac{x}{n_1}\right)\right|^2\frac{dx}{n_1}\\
&&\ll X^2\log^2  X
\sum_{\stackrel{m_1\leq X}{n_1\leq X}}
\frac{ |d_1(m_1,n_1)|X^{3/2}}{m_1 n_1^{1/4}}\ll X^{7/2}\log^2  X,
\end{eqnarray*}
which combining a splitting argument gives (4.21).

\section{Proof of Theorem 1}

In this section, we prove Theorem 1.
Define
\begin{eqnarray}\label{T-def}
T(m,n):=\sum_{\ell |gcd(m,n)}\ell \tau\left(\frac m\ell\right) \tau\left(\frac n\ell\right)
\end{eqnarray}
and
\begin{eqnarray}\label{S-def}
S(x):=\sum_{m,n\leq x}T(m,n).
\end{eqnarray}
The sum $S(x)$ plays an important role in the proof of Theorem 1.

\subsection{An asymptotic formula of $S(x)$}\

We easily see that
\begin{eqnarray}\label{S-tran}
S(x)&&=\sum_{\ell \max(j,k)\leq x}\ell \tau(j)\tau(k)\\
&&=\sum_{\stackrel{\ell k\leq x}{j\leq k}}\ell\tau(j)\tau(k)
+\sum_{\stackrel{\ell j\leq x}{k\leq j}}\ell\tau(j)\tau(k)
-\sum_{\ell k\leq x}\ell \tau^2(k)\nonumber\\
&&=2\sum_{\ell k\leq x}\ell\tau(k)\sum_{j\leq k}\tau(j)-\sum_{\ell k\leq x}\ell \tau^2(k).\nonumber
\end{eqnarray}

Inserting (\ref{1.1}) into (\ref{S-tran}), we get
\begin{eqnarray*}\label{SII-SI}
S(x)&&=2\sum_{\ell k\leq x}\ell\tau(k)\left(k\log k+(2\gamma-1)k+\Delta(k)\right)-\sum_{\ell k\leq x}\ell \tau^2(k)\\
&&=2\sum_{\ell k\leq x}\ell k\tau(k) \log k
 +(4\gamma-2)\sum_{\ell k\leq x}\ell k\tau(k)+
2\sum_{\ell k\leq x}\ell\tau(k) \Delta(k)-\sum_{\ell k\leq x}\ell \tau^2(k).\nonumber
\end{eqnarray*}

It is easy to see that
$$\frac 12\sum_{n=\ell k}\tau(k) \log k=
\frac 12\sum_{n=\ell n_1n_2}  (\log n_1+\log n_2)=\tau_3^{(1)}(n).
$$
Define
\begin{eqnarray}\label{s-j}
&&s_1(n):=n\tau_3^{(1)}(n), \ \ \ \ \ \ \ \mathcal{S}_1(x):=\sum_{n\leq x}s_1(n), \\
&& s_2(n):=n\tau_3(n),\ \ \ \ \ \ \ \ \ \ \mathcal{S}_2(x):=\sum_{n\leq x}s_2(n),\nonumber\\
&& s_3(n):= \sum_{n=\ell k}l\tau(k)\Delta(k),\ \ \ \
\mathcal{S}_3(x):=\sum_{n\leq x}s_3(n),\nonumber\\
&&  s_4(n):=\sum_{n=\ell k}\ell \tau^2(k),\ \ \ \ \ \ \ \ \ \mathcal{S}_4(x):=\sum_{n\leq x}s_4(n).\nonumber
 \end{eqnarray}

Let
\begin{equation}\label{s-def}
s(n):=4s_1(n)+(4\gamma-2)s_2(n)+2s_3(n)-s_4(n).
\end{equation}
Then we have
\begin{eqnarray}\label{Sx-def}
&&S(x)=\sum_{n\leq x}s(n)
=4 \mathcal{S}_1(x)+(4\gamma-2)\mathcal{S}_2(x)+2\mathcal{S}_3(x)-\mathcal{S}_4(x).
\end{eqnarray}

We first evaluate   $\mathcal{S}_3(x).$  By Lemma 2.8 and partial summation we have
\begin{eqnarray}\label{sum-of-t2}
 \mathcal{S}_3(x)&&= \sum_{ k\leq x}\tau(k) \Delta(k)\sum_{\ell\leq x/k}\ell\\
&&= \sum_{ k\leq x}\tau(k) \Delta(k)\left(\frac{x^2}{2k^2}+O(\frac xk)\right)\nonumber\\
&&=x^2\sum_{ k\leq x}\frac{\tau(k) \Delta(k)}{2k^2}+O\left(x\sum_{ k\leq x}\frac{|\tau(k)\Delta(k)|}{k }\right)\nonumber\\
&&=x^2\sum_{ k=1}^\infty\frac{\tau(k) \Delta(k)}{2k^2}+O(x^{5/4}\log^{5/2} x).\nonumber
\end{eqnarray}

Similarly by Lemma 2.6 with $\ell=2$ and partial summation we get
\begin{eqnarray}\label{sum-of-t3}
 \mathcal{S}_4(x)&&=\sum_{  k\leq x} \tau^2(k)\sum_{\ell\leq x/k}\ell\\
&&=\sum_{  k\leq x} \tau^2(k)\left(\frac{x^2}{2k^2}+O(\frac xk)\right)\nonumber\\
&&=\frac{x^2}{2}\sum_{  k=1}^\infty\frac{\tau^2(k)}{ k^2}+O(x\log^4 x).\nonumber
\end{eqnarray}

 By Lemma 2.10 we have
 \begin{eqnarray*}
D_3(x)=\sum_{n\leq x}\tau_3(n)=xC_2\log^2 x+C_1x\log x+C_0x+\Delta_3(x),
\end{eqnarray*}
with $\Delta_3(x)\ll x^{43/96+\varepsilon}.$
  So by partial integration we have
\begin{eqnarray}\label{ND3}
\ \ \mathcal{S}_2(x)&& =  \sum_{n\leq x}n\tau_3(n)=\int_{1^{-}}^x udD_3(u)\\
&&= \int_{1^{-}}^x ud(C_2u\log^2u+C_1u\log u+C_0u)+
\int_{1^{-}}^x u d\Delta_3(u) \nonumber\\
&&= \int_{1^{-}}^x u(C_2\log^2u+2C_2\log u+C_1\log u+C_1+C_0)du+O(x^{\frac{139}{96}+\varepsilon})\nonumber\\
&&=x^2(C_2^{*}\log^2 x+C_1^{*}\log x+C_0^{*})
+O(x^{\frac{139}{96}+\varepsilon}),\nonumber
\end{eqnarray}
with
$$ C_2^{*}= \frac{C_2}{2},\ \ C_1^{*}=\frac{C_1+C_2}{2},\ \   C_0^{*}=\frac{2C_0+C_1-C_2}{4}.$$

Similarly from (\ref{definition}), (\ref{residue}), Proposition 3.2  and partial integration we get
\begin{eqnarray}\label{ND31}
\mathcal{S}_1(x)&& =   \frac{K_3}{2}x^2\log^3 x+\left(\frac{2K_2+3K_3}{4}\right)x^2\log^2 x\\
&&\ \ \ +\left(\frac{2K_1+2K_2-3K_3}{4}\right)x^2\log x\nonumber\\
&&\ \ \ +
\left(\frac{4K_0+2K_1-2K_2+3K_3}{8}\right)x^2
+O(x^{\frac{139}{96}+\varepsilon}).\nonumber
\end{eqnarray}
where $K_j(j=0, 1, 2, 3)$ are constants defined in (\ref{residue}).

Collecting (\ref{SII-SI})-(\ref{ND31}) we get
\begin{eqnarray}\label{S-final}
S(x) =x^2\left(\sum_{j=0}^3d_j\log^j x\right)+O(x^{\frac{139}{96}+\varepsilon}),
\end{eqnarray}
where $d_j(j=0, 1, 2, 3)$ are computable constants.

\subsection{Proof of Theorem 1}\

We first consider $f(x)=S^{(1)}(x).$  From the relation(see Nowak and T\'{o}th \cite{NT})
$$s(m,n)=\sum_{d|gcd(m,n)}\mu(d)T(m/d,n/d)$$
we get
\begin{eqnarray}\label{S1-mid}
S^{(1)}(x)&&=
\sum_{m,n\leq x}s(m,n)=\sum_{d\leq x}\mu(d)\sum_{m,n\leq x/d}T(m,n)\\
&&=\sum_{d\leq x}\mu(d)S(x/d).\nonumber
\end{eqnarray}

Inserting (\ref{S-final}) into (\ref{S1-mid}) and by some easy calculations we get the asymptotic formula
\begin{eqnarray}\label{S1-final}
&&S^{(1)}(x)=x^2\left(\sum_{r=0}^3A_{r,S^{(1)}}\log^r x\right)+O(x^{\frac{139}{96}+\varepsilon}).
\end{eqnarray}

From the relation (see Nowak and T\'{o}th \cite{NT})
$$c(m,n)=\sum_{d|gcd(m,n)}\mu(d)s(m/d,n/d)$$
we get easily
\begin{eqnarray}\label{C1-mid}
C^{(1)}(x)&&=\sum_{m,n\leq x}c(m,n)=\sum_{d\leq x}\mu(d)\sum_{m,n\leq x/d}s(m,n)\\
&&=\sum_{d\leq x}\mu(d)S^{(1)}(x/d).\nonumber
\end{eqnarray}
Inserting (\ref{S1-final}) into (\ref{C1-mid}) and by some easy calculations we get the asymptotic formula
\begin{eqnarray}\label{C1-final}
&&C^{(1)}(x)=x^2\left(\sum_{r=0}^3A_{r,C^{(1)}}\log^r x\right)+O(x^{\frac{139}{96}+\varepsilon}).
\end{eqnarray}

\bigskip

Now we consider $S^{(2)}(x)$ and  $C^{(2)}(x). $ We have
\begin{eqnarray}\label{S-U}
&&S^{(2)}(x)=S^{(1)}(x)-U(x),\\
&&C^{(2)}(x)=C^{(1)}(x)-U(x),\nonumber
\end{eqnarray}
where $U(x)$ was defined by (\ref{4.1}). From (\ref{S1-final}), (\ref{C1-final}),
(\ref{S-U})  and
(\ref{4.3}) or Proposition 4.1,
  we get Theorem 1 for $S^{(2)}(x)$ and $C^{(2)}(x)$.

This completes the proof of Theorem 1.

\section{Proof of Theorem 2}\

In this section we prove Theorem 2.

\subsection{Properties of $R(s)$ }\

Define
\begin{eqnarray}\label{Rs-def}
R(s):=\sum_{n=1}^\infty\frac{\tau(n)\Delta(n)}{n^s}, \ \ \Re s>2.
\end{eqnarray}
In this subsection we study the function $R(s).$

For any $x\geq 1,$ define
$$Q(x):=\sum_{n\leq x}\tau(n)\Delta(n).$$
From Lemma 2.9, we can write
\begin{eqnarray}\label{6.2}
Q(x)=M_1(x)+E_Q(x),
\end{eqnarray}
where
\begin{eqnarray}\label{6.3}
&&M_1(x)=x\sum_{j=0}^3c_j(\log x)^j
\end{eqnarray}
and
\begin{eqnarray}\label{6.4}
&&\ \ E_Q(x)\ll x^{3/4}\log x.
\end{eqnarray}

From (\ref{6.2})-(\ref{6.4}) and partial summation it is easy to see  that the infinite series in the right-hand side of
 (\ref{Rs-def}) is absolutely convergent for $\Re s>5/4$ and we have the bound
  \begin{eqnarray}\label{6.5}
 R(s)\ll \frac{1}{5/4-\sigma},\ \ \ (\sigma>5/4).
 \end{eqnarray}

Without loss of generality, we suppose first $5/4<\Re s\leq 2$ and $|t|\geq 3.$
  Let $X\geq 2$ be a parameter to be chosen later. We write
 \begin{eqnarray}\label{6.6}
 R(s)=R_1(s;X)+R_2(s;X),
 \end{eqnarray}
where
\begin{equation}\label{6.7}
R_1(s;X):=\sum_{n\leq X}\frac{\tau(n)\Delta(n)}{n^s},\ \ \
R_2(s;X):=\sum_{n> X}\frac{\tau(n)\Delta(n)}{n^s}
\end{equation}

It is easy to see that  $R_1(s;X)$ is an entire function on ${\Bbb C}.$

 By partial summation we have that
 \begin{eqnarray}\label{6.8}
 R_2(s;X)=\int_X^\infty \frac{1}{u^s}dQ(u)=
 \int_X^\infty \frac{1}{u^s}dM_1(u)+\int_X^\infty \frac{1}{u^s}dE_Q(u).
\end{eqnarray}
 From  (\ref{6.3}) it is easy to check that
 $$M_1^{\prime}(u)=\sum_{j=0}^3c_j^{*}\log^j u$$
 with $c_3^{*}=c_3$ and $  c_j^{*}=c_j+c_{j+1}(j+1)\ \ (j=0,1,2).$

 So by partial integration and easy calculations we have
 \begin{eqnarray}  \label{6.9}
 &&\ \ \int_X^\infty \frac{1}{u^s}dM_1(u) =\int_X^\infty \frac{1}{u^s}M_1^{\prime}(u)du\\
 &&=\frac{6c_3^{*}X^{1-s} }{(s-1)^4}+\frac{X^{1-s}}{(s-1)^3}
 \left(6c_3^{*}\log X+2c_2^{*}\right)\nonumber\\
 &&\ \ +\frac{X^{1-s}}{(s-1)^2}\left(3c_3^{*}\log^2 X+2c_2^{*}\log X
 +c_1^{*}\right)\nonumber\\
 &&\ \ + \frac{X^{1-s}}{s-1}
 \left(c_3^{*}\log^3 X +c_2^{*}\log^2 X  +c_1^{*}\log X+  c_0^{*}  \right),\nonumber
\end{eqnarray}
which is meromorphic  for $s\in {\Bbb C}$ and
  $s=1$ is a pole of order $4.$

   By partial integration we get
 \begin{eqnarray}\label{6.10}
\int_X^\infty \frac{1}{u^s}dE_Q(u)=\frac{E_Q(X)}{X^s}+s\int_X^\infty\frac{E_Q(u)}{u^{s+1}}du.
\end{eqnarray}
The estimate (\ref{6.4}) implies that the integral in the right-hand side of
(\ref{6.10}) is absolutely convergent for $\sigma>3/4.$ From (\ref{6.6})-(\ref{6.10}) we see that $R(s)$
can be continued meromorphically to $\Re s>3/4,$ and $s=1$ is a pole of order $4.$

Suppose now $3/4<\sigma\leq 5/4$ and $|t|\geq 3.$ From Lemma 2.8 and partial summation we have
\begin{eqnarray} \label{6.11}
 R_1(s;X)\ll X^{5/4-\sigma} \log X.
\end{eqnarray}
 From (\ref{6.9}) we get
 \begin{eqnarray}\label{6.12}
 \int_X^\infty \frac{1}{u^s}dM_1(u)\ll \frac{X^{1-\sigma}}{|t|} \log^3 X.
\end{eqnarray}
 From (\ref{6.4}) and (\ref{6.10}) we have
  \begin{eqnarray}\label{6.13}
\int_X^\infty \frac{1}{u^s}dE_Q(u)\ll X^{3/4-\sigma}\log X+|t|X^{3/4-\sigma}\ll |t|X^{3/4-\sigma}.
\end{eqnarray}

 From (\ref{6.6}), (\ref{6.8}),   (\ref{6.11})-(\ref{6.13}) and   choosing $X=t^2$ we get that
 \begin{eqnarray}\label{Rs-bound}
 R(\sigma+it)\ll |t|^{5/2-2\sigma}\log |t|,\ \ 3/4<\sigma\leq 5/4.
 \end{eqnarray}

\subsection{Proof of Theorem 2}\
We first prove Theorem 2 for $S^{(1)}(x)$ and $C^{(1)}(x).$

Define
\begin{equation*}
\mu_1(n):=\mu(n),\ \ \ \ \mu_2(n):=\sum_{n=n_1n_2}\mu(n_1)\mu(n_2).
\end{equation*}

For $j=1,2,$ define
\begin{eqnarray}\label{t1j-def}
&&t_{1,j}(n):=\sum_{n=dm}\mu_j(d)s_1(m),\ \ \ \ \
\mathcal{T}_{1,j}(x)=\sum^{\prime}_{n\leq x}t_{1,j}(n),\\
&&t_{2,j}(n):=\sum_{n=dm}\mu_j(d)s_2(m),\ \ \ \ \ \
\mathcal{T}_{2,j}(x)=\sum_{n\leq x}^{\prime}t_{2,j}(n), \nonumber\\
&&t_{3,j}(n):=\sum_{n=dm}\mu_j(d)s_3(m),  \ \ \ \ \ \ \ \  \mathcal{T}_{3,j}(x)=\sum_{n\leq x}^{\prime}t_{3,j}(n),  \nonumber\\ &&
t_{4,j}(n):=\sum_{n=dm}\mu_j(d)s_4(m),\ \ \ \ \ \ \  \ \mathcal{T}_{4,j}(x)=\sum_{n\leq x}^{\prime}t_{4,j}(n), \nonumber
\end{eqnarray}
where $\sum^{\prime}_{n\leq x}f(n)$ means that if $n=x$ is an integer, then
this term is counted with only $f(x)/2$, and where $s_j(n)(j=1, 2, 3, 4)$ were
defined in (\ref{s-j}). Define
\begin{equation}\label{t-def}
t_j(n):=\sum_{n=dm}\mu_j(d)s(m),\ \ (j=1,2)
\end{equation}
where  $s(m)$ was defined by (\ref{s-def}).

From (\ref{s-def}) and (\ref{t-def}), we have
\begin{equation}\label{tj-def}
t_j(n)=4t_{1,j}(n)+(4\gamma-2)t_{2,j}(n)+2t_{3,j}(n)-t_{4,j}(n),\ \ j=1,2.
\end{equation}

Define further $(j=1,2)$
\begin{eqnarray}\label{Fj-def}
&&\mathfrak{F}_{1,j}(s):=\sum_{n=1}^\infty\frac{t_{1,j}(n)}{n^s},   \ \ \ \
 \mathfrak{F}_{2,j}(s):=\sum_{n=1}^\infty\frac{t_{2,j}(n)}{n^s},  \ \ \
 \mathfrak{F}_{3,j}(s):=\sum_{n=1}^\infty\frac{t_{3,j}(n)}{n^s}, \\&&
\mathfrak{F}_{4,j}(s):=\sum_{n=1}^\infty\frac{t_{4,j}(n)}{n^s},\ \ \
\mathfrak{F}_j(s):=\sum_{n=1}^\infty\frac{t_{j}(n)}{n^s}. \nonumber
\end{eqnarray}

From (\ref{tj-def}) and (\ref{Fj-def}) we have
\begin{equation}\label{Fj-s}
\mathfrak{F}_j(s)=4\mathfrak{F}_{1,j}(s)+(4\gamma-2)\mathfrak{F}_{2,j}(s)
+2\mathfrak{F}_{3,j}(s)-\mathfrak{F}_{4,j}(s),\ \ j=1,2.
\end{equation}

From Subsection 5.2 we see that ($j=1,2$)
\begin{eqnarray}\label{D-series}
&&\mathfrak{F}_{1,j}(s)=-\frac{\zeta^2(s-1)\zeta^{\prime}(s-1)}{\zeta^j(s)},   \ \ \ \
 \mathfrak{F}_{2,j}(s)=\frac{\zeta^3(s-1) }{\zeta^j(s)} ,  \\
 &&\mathfrak{F}_{3,j}(s)=\frac{\zeta(s-1)R(s) }{\zeta^j(s)} ,  \ \ \ \ \
\mathfrak{F}_{4,j}(s)= \frac{\zeta(s-1)\zeta^{4-j}(s) }{\zeta(2s)}, \nonumber
\end{eqnarray}
where $R(s)$ was defined in (\ref{Rs-def}), and in the expression of $\mathfrak{F}_{4,j}(s)$
we used the fact (see, for example, Chapter 1 of Ivi\'c \cite{Iv})
$$\sum_{n=1}^\infty\frac{\tau^2(n)}{n^s}=\frac{\zeta^4(s)}{\zeta(2s)}\ \ (\Re s>1).$$

Define further
\begin{eqnarray}\label{2-residue}
&&\mathfrak{M}_{1,j}(x):=Res_{s=2}\mathfrak{F}_{1,j}(s)x^ss^{-1},
  \ \ \ \ \
 \mathfrak{M}_{2,j}(x):=Res_{s=2}\mathfrak{F}_{2,j}(s)x^ss^{-1}, \\
 &&   \mathfrak{M}_{3,j}(x):=Res_{s=2}\mathfrak{F}_{3,j}(s)x^ss^{-1},
  \ \ \ \ \
\mathfrak{M}_{4,j}(x):=Res_{s=2}\mathfrak{F}_{4,j}(s)x^ss^{-1}.\nonumber
\end{eqnarray}

 Finally we define
 \begin{eqnarray} \label{error-j}
\Delta_{\mathcal{T}_{i,j}}(x):=\mathcal{T}_{i,j}(x)-\mathfrak{M}_{i,j}(x),\ \ \ 1\leq i\leq 4, j=1,2.
\end{eqnarray}

From (\ref{t1j-def})-(\ref{Fj-s}), (\ref{2-residue}), (\ref{error-j}) and (\ref{S1-mid})  we see that
$$ 4\mathfrak{M}_{1,1}(x) +(4\gamma-2)\mathfrak{M}_{2,1}(x) +2\mathfrak{M}_{3,1}(x) -\mathfrak{M}_{4,1}(x)  $$ is the main term in the asymptotic formula of $S^{(1)}(x)$
and we have
\begin{eqnarray*}\label{S1-error}
\ \ \ \ \ \ E_{S^{(1)}}(x)&&=4\Delta_{ \mathcal{T}_{1,1}}(x)
+(4\gamma-2)\Delta_{\mathcal{T}_{2,1}}(x)
+2\Delta_{ \mathcal{T}_{3,1}}(x)-\Delta_{ \mathcal{T}_{4,1}}(x)+O(x^{1+\varepsilon}).
\end{eqnarray*}
Note that  the term $O(x^{1+\varepsilon})$ appears because it contains
$O(|t_{1,1}^{(1)}(x)|+|t_{2,2}(x)|+|t_{3,2}(x)|+|t_{4,2}(x)|)$ when $x$ is an integer.
From (\ref{sum-of-t3}) and the definition of $t_{4,1}(n),$ it is easy to check that
$\Delta_{ \mathcal{T}_{4,1}}(x)\ll x\log^4 x.$ Thus we have
\begin{eqnarray}\label{S1-error}
\ \ \ \ \ \ E_{S^{(1)}}(x)&&=4\Delta_{T_{1,1}}(x)+(4\gamma-2)\Delta_{T_{2,1}}(x)
+2\Delta_{T_{3,1}}(x) +O(x^{1+\varepsilon}).
\end{eqnarray}

Similarly, $ 4\mathfrak{M}_{1,2}(x) +(4\gamma-2)\mathfrak{M}_{2,2}(x) +2\mathfrak{M}_{3,2}(x) -\mathfrak{M}_{4,2}(x)  $ is the main term in the asymptotic formula of $C^{(1)}(x)$
and we have
\begin{equation}\label{C1-error}
E_{C^{(1)}}(x)=4\Delta_{\mathcal{T}_{1,2}}(x)+(4\gamma-2)\Delta_{\mathcal{T}_{2,2}}(x)
+2\Delta_{\mathcal{T}_{3,2}}(x) +O(x^{1+\varepsilon}).
\end{equation}

By Perron's formula    we easily see that
\begin{eqnarray}\label{M1m1}
\int_0^TE_{S^{(1)}}(x)dx=\frac{1}{2\pi i}\int_{2-\varepsilon-\infty}^{2-\varepsilon+\infty}
\frac{\mathfrak{F}_1(s)T^{1+s}}{s(s+1)}ds.
\end{eqnarray}

By Lemma 2.3 and  (\ref{Rs-bound}) it is easy to see that
\begin{equation}\label{Fj-76}
 \mathfrak{F}_j(7/6+\varepsilon+it)\ll |t|^{1-\varepsilon},\ \ |t|\geq 2.
 \end{equation}

Moving the integration line to
$\Re s=7/6+\varepsilon$ in the right-hand side of (\ref{M1m1}) and then using (\ref{Fj-76}) we get
\begin{eqnarray}\label{S1-mean}
\int_0^TE_{S^{(1)}}(x)dx&&\ll T^{13/6+\varepsilon} \int_{7/6+\varepsilon-\infty}^{7/6+\varepsilon+\infty}
\frac{|\mathfrak{F}_1(7/6+\varepsilon+it)| }{|t|^2+1}dt\\
&&\ll T^{13/6+\varepsilon}.\nonumber
\end{eqnarray}
Similarly
\begin{eqnarray}\label{C1-mean}
\int_0^TE_{C^{(1)}}(x)dx \ll   T^{13/6+\varepsilon}.
\end{eqnarray}

 Now we consider $E_{S^{(2)}}(x)$ and $E_{C^{(2)}}(x)$.  From
 (\ref{S-U}), (\ref{S1-error}) and (\ref{C1-error}) we get that
\begin{eqnarray}\label{S2-error-expr}
&&E_{S^{(2)}}(x)=4\Delta_{\mathcal{T}_{1,1}}(x)+(4\gamma-2)\Delta_{\mathcal{T}_{2,1}}(x)
+2\Delta_{\mathcal{T}_{3,1}}(x)-E_U(x) +O(x^{1+\varepsilon}),\\
 &&
E_{C^{(2)}}(x)=4\Delta_{\mathcal{T}_{1,2}}(x)+(4\gamma-2)\Delta_{\mathcal{T}_{2,2}}(x)
+2\Delta_{\mathcal{T}_{3,2}}(x)-E_U(x) +O(x^{1+\varepsilon}).
\end{eqnarray}

From (\ref{S1-mean}), (\ref{C1-mean}) and Proposition 4.2 we get
\begin{eqnarray*}
\int_0^TE_{S^{(2)}}(x)dx \ll T^{13/6+\varepsilon},\ \ \ \
\int_0^TE_{C^{(2)}}(x)dx \ll T^{13/6+\varepsilon}.
\end{eqnarray*}

This completes the proof of Theorem 2.

\section{Some moment results involving $\zeta(s)$}\
Throughout this section, suppose $7/12<\sigma<1$ is a fixed real number and
 $T\geq 3$ is a large parameter. In this  section, we study some
 integrals involving $\zeta(s),$ which are important to the proof of Theorem 3.

For any fixed integer $j\geq 1,$ define the following integrals
\begin{eqnarray}\label{integrals}
&&I_{1,j}(T;\sigma):=\int_1^{T}
\left|\frac{\zeta^2(\sigma+it)\zeta^{\prime}(\sigma+it)}
{\zeta^j(2- \sigma+it)}\right|^2 dt,\\
&&I_{2,j}(T;\sigma):=\int_1^{T}
\left|\frac{\zeta^3(\sigma+it)}{\zeta^j(2- \sigma+it)} \right|^{2}dt,\nonumber\\
&&I_{3,j}(T;\sigma):=\int_1^T \frac{\zeta^2(\sigma-it)\zeta^{\prime}(\sigma-it)}
{\zeta^j(1+ \sigma-it)}\times
\frac{\zeta^3(\sigma+it)}{\zeta^j(1+ \sigma+it)}
dt.\nonumber
\end{eqnarray}

\subsection{Evaluations of $I_{1,j}(T;\sigma)$ and $I_{2,j}(T;\sigma)$}\

We consider $I_{1,j}(T;\sigma)$ first.
From the trivial relation (recall (\ref{diri-series}))
$$ (\zeta^2(s)\zeta^{\prime}(s))^2=
 \left(\sum_{n\leq T}\tau_3^{(1)}(n)n^{-s} \right)^2+
(\zeta^2(s)\zeta^{\prime}(s))^2 - \left(\sum_{n\leq T} \tau_3^{(1)}(n)n^{-s}\right)^2
$$

we have
\begin{eqnarray*}
&&\left(\frac{\zeta^2(\sigma+it )\zeta^{\prime}(\sigma+it )}
{\zeta^{j}(2- \sigma+it)}\right)^2=\frac{\left(\sum_{n\leq T}\tau_3^{(1)}(n)n^{-\sigma-it}\right)^2}
{\zeta^{2j}(2- \sigma+it )}\\
&&\ \ \ +\frac{1}{\zeta^{2j}(2- \sigma+it)}\times
\left((\zeta^2(\sigma+it )\zeta^{\prime}(\sigma+it))^2 - \left(\sum_{n\leq T} \tau_3^{(1)}(n)n^{-\sigma-it }\right)^2  \right).
\end{eqnarray*}
Thus we get that (note that $1/\zeta(2-\sigma+it)\ll 1$)
\begin{eqnarray}\label{tran}
&&\ \ \ \ \ \ \ I_{1,j}(T;\sigma)=J_{11}(T;\sigma)+ O(J_{12}(T;\sigma)),\\
&& J_{11}(T;\sigma):=\int_1^{T}\left|\frac{\sum_{n\leq T}\tau_3^{(1)}(n)n^{-\sigma-it}}{\zeta^j(2- \sigma+it)}\right|^2dt,\nonumber\\
&&J_{12}(T;\sigma):=\int_1^T\left|(\zeta^2(\sigma+it )\zeta^{\prime}(\sigma+it))^2 - \left(\sum_{n\leq T} \tau_3^{(1)}(n)n^{-\sigma-it }\right)^2 \right|dt.\nonumber
\end{eqnarray}

We first bound the  integral $J_{12}(T;\sigma).$    Let now
$$F(s)=  (\zeta^2(s)\zeta^{\prime}(s) )^2
 -\left(\sum_{n\leq T}\tau_3^{(1)}(n)n^{-s}\right)^2,$$
and apply Lemma 2.12 with $q=1, \alpha=7/12+\delta, \beta=1+\delta, \gamma=\sigma,$
where $0<\delta<1/2$ is a fixed constant which may be chosen arbitrarily small.

By Lemma 2.4 with $j=0,1$ and H\"{o}lder's inequality we get
\begin{eqnarray*}
\int_1^T\left| \zeta^2(\alpha+it)\zeta^{\prime}(\alpha+it) \right|^2dt\ll T^{1+\delta}.
\end{eqnarray*}
From Lemma 2.11 we have
\begin{eqnarray*}
\int_1^T\left| \sum_{n\leq T}\tau_3^{(1)}(n)n^{-\alpha-it}\right|^2dt\ll T.
\end{eqnarray*}

From the above two estimates we get
\begin{equation}\label{alpha-bound}
\int_1^T |F(\alpha+it)|dt\ll T^{1+\delta}.
\end{equation}

For $F(\beta+it),$ it is easy to see that
\begin{eqnarray*}
F(\beta+it)&&=\left( \zeta^2(1+\delta+it)\zeta^{\prime}(1+\delta+it) \right)^2-\left(\sum_{n\leq T}\tau_3^{(1)}(n)n^{-1-\delta-it}\right)^2\\
&&=\sum_{n>T}g(n,T)n^{-1-\delta-it},
\end{eqnarray*}
where $|g(n,T)|\leq a_3^{(1)}*a_3^{(1)}(n)\leq \tau_6(n)\log^2 n.$ Using Lemma 2.11 and Cauchy's inequality we get
\begin{equation}\label{beta-bound}
\int_1^T |F(\beta+it)|dt\ll T^{1/2}\left(\int_1^T |\sum_{n>T}g(n,T)n^{-1-\delta-it}|^2dt\right)^{1/2}\ll T^{1/2}.
\end{equation}

From (\ref{alpha-bound}), (\ref{beta-bound}) and Lemma 2.12 we get for
$7/12<\sigma<1$
\begin{eqnarray}\label{err-r-bound}
&&\ \ \ J_{12}(T;\sigma)= \int_1^T |F(\sigma+it)|dt\ll T^A,\\
A&&=(1+\delta)\frac{1+\delta-\sigma}{\frac{5}{12}}+
\frac 12\times \frac{\sigma-7/12-\delta}{\frac{5}{12}}\nonumber\\
&&\ \leq \frac{17-12\sigma}{10}+3\delta\leq \frac{17-12\sigma}{10}+\varepsilon\nonumber
\end{eqnarray}
if $\delta\leq \varepsilon/3.$

Now we evaluate the first integral $J_{11}(T;\sigma)$ in the right-hand side of (\ref{tran}). Let
$\mu_j(\ell)$ be defined by
$$\sum_{\ell=1}^\infty\frac{\mu_j(\ell)}{\ell^{z}}=\frac{1}{\zeta^j(z)},\ \ \Re z>1.$$
Obviously
$$\mu_1(\ell)=\mu(\ell),\ \ \mu_j(\ell)=\sum_{\ell=\ell_1\cdots \ell_j}\mu(\ell_1)\cdots \mu(\ell_j)\ (j\geq 2). $$
Write
\begin{eqnarray}\label{1-expression}
\frac{\sum_{n\leq T}\tau_3^{(1)}(n)n^{-\sigma-it}}{\zeta^j(2- \sigma+it)}
=\sum_{n\leq T}\frac{\tau_3^{(1)}(n)}{n^{\sigma+it}}
\sum_{\ell=1}^\infty\frac{\mu_j(\ell)}{\ell^{2-\sigma+it}}=
\sum_{m=1}^\infty \frac{a_{3,j}^{(1)}(m;\sigma,T)}{m^{\sigma+it}},
\end{eqnarray}
where
\begin{eqnarray}\label{function-def}
a_{3,j}^{(1)}(m;\sigma,T):=\sum_{\stackrel{m=\ell n}{n\leq T}}\frac{\tau_3^{(1)}(n)\mu_j(\ell)}{\ell^{2-2\sigma}},\ \
a_{3,j}^{(1)}(m;\sigma):=\sum_{ m=\ell n}\frac{\tau_3^{(1)}(n)\mu_j(\ell)}{\ell^{2-2\sigma}}.
\end{eqnarray}
Obviously we have
$a_{3,j}^{(1)}(m;\sigma,T)=a_{3,j}^{(1)}(m;\sigma),\ \ \ m\leq T.$

By Lemma 2.11 we have
\begin{eqnarray}\label{111-mid}
J_{11}(T;\sigma)
=T\sum_{m=1}^\infty\frac{(a_{3,j}^{(1)}(m;\sigma,T))^2}{m^{2\sigma}}
  +O\left(\sum_{m=1}^\infty\frac{(a_{3,j}^{(1)}(m;\sigma,T))^2}
{m^{2\sigma-1}}\right).
\end{eqnarray}

By the simple inequalities
\begin{equation}\label{inequality}
|a_{3,j}^{(1)}(m;\sigma,T)|\leq \tau_{3+j}(m)\log m,
\ \ |a_{3,j}^{(1)}(m;\sigma)|\leq \tau_{3+j}(m)\log m,
\end{equation}
 we get
with the help of Lemma 2.6 by taking $\ell=3+j$ that
\begin{eqnarray}\label{main-term-cons}
\ \ \sum_{m=1}^\infty\frac{(a_{3,j}^{(1)}(m;\sigma,T))^2}{m^{2\sigma}}
&&=\sum_{m=1}^\infty\frac{(a_{3,j}^{(1)}(m;\sigma))^2}{m^{2\sigma}}
+O\left(\sum_{m>T}\frac{\tau_{3+j}^2(m)\log^2 m }{m^{2\sigma}}\right)\\
&&=\sum_{m=1}^\infty\frac{(a_{3,j}^{(1)}(m;\sigma))^2}{m^{2\sigma}}
+O(T^{2-2\sigma}(\log T)^{(3+j)^2+1}).\nonumber
\end{eqnarray}

Write
$$\sum_{m=1}^\infty\frac{(a_{3,j}^{(1)}(m;\sigma,T))^2}
{m^{2\sigma-1}}=\sum_{m\leq T}\frac{(a_{3,j}^{(1)}(m;\sigma,T))^2}
{m^{2\sigma-1}}+\sum_{m>T} \frac{(a_{3,j}^{(1)}(m;\sigma,T))^2}
{m^{2\sigma-1}}.$$

By (\ref{inequality}) and Lemma 2.6 with $\ell=3+j$ again we have
\begin{equation}\label{former}
\sum_{m\leq T}\frac{(a_{3,j}^{(1)}(m;\sigma,T))^2}
{m^{2\sigma-1}}\ll T^{2-2\sigma}(\log T)^{(3+j)^2+1}.
\end{equation}

From the first definition in (\ref{function-def})  we have
\begin{eqnarray*} \label{division222}
&&\ \ \ \ \ \sum_{m>T} \frac{(a_{3,j}^{(1)}(m;\sigma,T))^2}
{m^{2\sigma-1}}\\&&\ll \sum_{\stackrel{\ell_1 n_2=\ell_2 n_1>T}{n_1, n_2\leq T}}
\frac{|\mu_j(\ell_1)\mu_j(\ell_2)|\tau_3^{(1)}(n_1)\tau_3^{(1)}(n_2)}
{\ell_1^{2-2\sigma}\ell_2^{2-2\sigma} }\times \frac{1}{(\ell_1 n_2)^{\sigma-\frac 12}   (\ell_2 n_1)^{\sigma-\frac 12} }\nonumber\\
&&=\sum_{\stackrel{\ell_1 n_2=\ell_2 n_1>T}{n_1, n_2\leq T}}
\frac{|\mu_j(\ell_1)\mu_j(\ell_2)|\tau_3^{(1)}(n_1)\tau_3^{(1)}(n_2)}
{\ell_1^{\frac 32- \sigma}\ell_2^{\frac 32-\sigma} n_1^{\sigma-\frac 12}
n_2^{\sigma-\frac 12}} \nonumber.
 \end{eqnarray*}

The equality $\ell_1/n_1=\ell_2/n_2\ (\ell_1, \ell_2, n_1, n_2\in {\Bbb N})$ implies that there exist
$\delta\in {\Bbb N}$ and $  \rho\in {\Bbb N}$ with $gcd(\delta, \rho)=1,$
such that
\begin{equation*}
\ell_1=\delta t,\ n_1=\rho t,\ \ell_2=\delta v,\ n_2=\rho v,\ t,v\in {\Bbb N}.
\end{equation*}

So we have
\begin{eqnarray} \label{division222}
\ \ \ \ \ \sum_{m>T} \frac{(a_{3,j}^{(1)}(m;\sigma,T))^2}
{m^{2\sigma-1}} &&
\ll \sum_{\stackrel{\delta tv \rho>T}{\rho\leq T/t, \rho\leq T/v}}
\frac{|\mu_j(\delta t)\mu_j( \delta v)|\tau_3^{(1)}(\rho t)\tau_3^{(1)}( \rho v)}
{\delta^{3-2\sigma}tv \rho^{2\sigma-1}}  \\
&&\ll \sum_{ \delta, t, v  }
\frac{|\mu_j(\delta t)\mu_j( \delta v)|\tau_3^{(1)}(t)\tau_3^{(1)}(v)(1+\log t)(1+\log v)}
{\delta^{3-2\sigma}tv}\nonumber\\
 &&\ \ \ \ \ \ \ \ \ \ \times \sum_{\rho\leq \min\left(\frac{T}{t},\frac{T}{v}\right)}
 \frac{\tau_3^2(\rho)(\log \rho+1)^2}{\rho^{2\sigma-1}},\nonumber
 \end{eqnarray}
 where we used the estimate
 $$\tau_3^{(1)}(\rho h)\leq \tau_3(\rho h) \log \rho h\leq
 \tau_3(\rho)\tau_3(h) (\log \rho+1)(\log  h+1), \ \ h\in {\Bbb N}.
 $$

 By Lemma 2.6 with $\ell=3$ we see that the innermost sum in the last line of
 (\ref{division222})   is
 \begin{eqnarray}\label{innermost}
\sum_{\rho\leq \min\left(\frac{T}{t},\frac{T}{v}\right)}
 \frac{\tau_3^2(\rho)(\log \rho+1)^2}{\rho^{2\sigma-1}}&&\ll
 \min\left(\frac{T}{t},\frac{T}{v}\right)^{2-2\sigma}(\log T)^{10}\\
 &&\ll \frac{T^{2-2\sigma}}{t^{1-\sigma}v^{1-\sigma}}(\log T)^{10},\nonumber
 \end{eqnarray}
 where we used the inequality $\min(x,y)\leq x^{1/2}y^{1/2},\ (x>0, y>0).$

 From (\ref{division222})  and (\ref{innermost})  we get
 \begin{eqnarray}\label{latter}
&&\ \ \ \ \ \ \ \sum_{m>T} \frac{(a_{3,j}^{(1)}(m;\sigma,T))^2}
{m^{2\sigma-1}}\\&&\ll T^{2-2\sigma}(\log T)^{10}
\sum_{ \delta, t, v  }
\frac{|\mu_j(\delta t)\mu_j( \delta v)|\tau_3^{(1)}(t)\tau_3^{(1)}(v)(1+\log t)(1+\log v)}{\delta^{3-2\sigma}t^{2-\sigma}v^{2-\sigma}}\nonumber\\
&&\ll T^{2-2\sigma}(\log T)^{10} \nonumber
\end{eqnarray}
 by noting that $7/12<\sigma<1.$

So (\ref{former}) and (\ref{latter}) we get the estimate
\begin{eqnarray}\label{second-series}
 \sum_{m=1}^\infty \frac{(a_{3,j}^{(1)}(m;\sigma,T))^2}
{m^{2\sigma-1}} \ll T^{2-2\sigma}(\log T)^{(3+j)^2+1}.
\end{eqnarray}

Inserting (\ref{main-term-cons}) and (\ref{second-series}) into (\ref{111-mid}), we get
\begin{eqnarray}\label{111-middle}
 J_{11}(T;\sigma)= B_{1,j}(\sigma)T+O(T^{2-2\sigma}(\log T)^{(3+j)^2+1}),
\end{eqnarray}
 where
 \begin{equation}
B_{1,j}(\sigma):= \sum_{m=1}^\infty\frac{(a_{3,j}^{(1)}(m;\sigma))^2}{m^{2\sigma}}.
\end{equation}

From (\ref{err-r-bound}) and (\ref{111-middle}) we finally get the asymptotic formula
\begin{eqnarray}\label{111-final}
I_{1,j}(T;\sigma) =B_{1,j}(\sigma)T+O(T^{(17-12\sigma)/10+\varepsilon}),\ \ \ \ (7/12<\sigma<1).
\end{eqnarray}

Similarly, we can prove that
\begin{eqnarray}\label{222-final}
I_{2,j}(T;\sigma) =B_{2,j}(\sigma)T+O(T^{(17-12\sigma)/10+\varepsilon}),\ \ \ \ (7/12<\sigma<1)
\end{eqnarray}
where
 \begin{equation}\label{I2j-cons}
B_{2,j}(\sigma):= \sum_{m=1}^\infty\frac{(a_{3,j}(m;\sigma))^2}{m^{2\sigma}},
\ \ \ \ a_{3,j}(m;\sigma):=\sum_{m=\ell n}\frac{\tau_3(n)\mu_j(\ell)}{\ell^{2-2\sigma}}.
\end{equation}
 We omit the proof of (\ref{222-final}) since it is almost the same as that of (\ref{111-final}).

 \subsection{Evaluation of $I_{3,j}(T;\sigma)$}\

Now we study $I_{3,j}(T;\sigma).$ From the decomposition (recall (\ref{diri-series}))
\begin{eqnarray*}
 &&\zeta^2(\sigma-it)\zeta^{\prime}(\sigma-it)
 \times
 \zeta^3(\sigma+it)\\
 && =
-\sum_{m\leq T}  \frac{\tau_3^{(1)}(m)}{m^{\sigma-it}}
\sum_{n\leq T} \frac{\tau_3(n)}{n^{\sigma+it}}\\
&&\ \ +\left( \zeta^2(\sigma-it)\zeta^{\prime}(\sigma-it)
 +\sum_{m\leq T}  \frac{\tau_3^{(1)}(m)}{m^{\sigma-it}}\right)
 \zeta^3(\sigma+it) \nonumber\\
&&\ \  +\left(-\sum_{m\leq T}  \frac{\tau_3^{(1)}(m)}{m^{\sigma-it}}\right)
\left( \zeta^3(\sigma+it) -\sum_{n\leq T} \frac{\tau_3(n)}{n^{\sigma+it}}
\right)
\end{eqnarray*}
we get  that
\begin{eqnarray*}
&&\ \ \ \ \   \frac{\zeta^2(\sigma-it)\zeta^{\prime}(\sigma-it)}
{\zeta^j(2- \sigma-it)} \times
\frac{\zeta^3(\sigma+it)}{\zeta^j(2- \sigma+it)}\\
&&=\frac{ -\sum_{m\leq T}  \frac{\tau_3^{(1)}(m)}{m^{\sigma-it}} }
{\zeta^j(2- \sigma-it)} \times\frac{
 \sum_{n\leq T} \frac{\tau_3(n)}{n^{\sigma+it}} }
{\zeta^j(2- \sigma+it)}\nonumber\\
&&\ \ \ +\left(\zeta^2(\sigma-it)\zeta^{\prime}(\sigma-it)
 +\sum_{m\leq T}  \frac{\tau_3^{(1)}(m)}{m^{\sigma-it}} \right)\times
\frac{\zeta^3(\sigma+it)}{|\zeta(2- \sigma+it)|^{2j}}\nonumber\\
&&\ \ \ +\left(-\sum_{m\leq T}  \frac{a_3^{(1)}(m)}{m^{\sigma-it}}\right)\times
\left( \zeta^3(\sigma+it) -\sum_{n\leq T} \frac{\tau_3(n)}{n^{\sigma+it}}\right)\times \frac{1}{|\zeta(2-\sigma+it)|^{2j}}.\nonumber
\end{eqnarray*}

Thus we have
\begin{eqnarray}\label{T3j-division}
I_{3,j}(T;\sigma)=J_{31}(T,\sigma)+O(J_{32}(T,\sigma)+J_{33}(T,\sigma) ),
\end{eqnarray}
where
\begin{eqnarray*}
&&J_{31}(T,\sigma): =\int_1^T \frac{ -\sum_{m\leq T}  \frac{\tau_3^{(1)}(m)}{m^{\sigma-it}} }{\zeta^j(2- \sigma-it)} \times\frac{ \sum_{n\leq T} \frac{\tau_3(n)}{n^{\sigma+it}} }
{\zeta^j(2- \sigma+it)}dt,\\
&& J_{32}(T,\sigma): =\int_1^T \left|\zeta^2(\sigma-it)\zeta^{\prime}(\sigma-it)
 +\sum_{m\leq T}  \frac{\tau_3^{(1)}(m)}{m^{\sigma-it}} \right|\times
|\zeta^3(\sigma+it)|  dt,\\
&&J_{33}(T,\sigma): =\int_1^T  \left|\sum_{m\leq T}  \frac{\tau_3^{(1)}(m)}{m^{\sigma-it}}\right|\times
\left| \zeta^3(\sigma+it) -\sum_{n\leq T} \frac{\tau_3(n)}{n^{\sigma+it}}\right|dt.
\end{eqnarray*}

We first consider $J_{31}(T;\sigma).$ Define
\begin{equation}\label{tau-3-j-function}
a_{3,j}(m;\sigma,T):=\sum_{\stackrel{m=\ell n}{n\leq T}}\frac{\tau_3(n)\mu_j(\ell)}{\ell^{2-2\sigma}}.
\end{equation}
Then we can write
\begin{eqnarray*}
&&\ \ \ \ \ \frac{ -\sum_{m\leq T}  \frac{\tau_3^{(1)}(m)}{m^{\sigma-it}} }
{\zeta^j(2- \sigma-it)} \times\frac{ \sum_{n\leq T} \frac{\tau_3(n)}{n^{\sigma+it}} }
{\zeta^j(2- \sigma+it)}\\
&&=\sum_{m_1=1}^\infty\frac{-a_{3,j}^{(1)}(m_1;\sigma,T)}{m_1^{\sigma}}
\sum_{m_2=1}^\infty\frac{a_{3,j}(m_2;\sigma,T)}
{m_2^{\sigma}}\left(\frac{m_1}{m_2}\right)^{it}\\
&&=\sum_{m=1}^\infty\frac{-a_{3,j}^{(1)}(m;\sigma,T)a_{3,j}(m;\sigma,T) }{m^{2\sigma}}
 +\sum_{\stackrel{m_1,m_2\in {\Bbb N}}{m_1\not= m_2}} \frac{-a_{3,j}^{(1)}(m_1;\sigma,T)a_{3,j}(m_2;\sigma,T)}
 {m_1^{\sigma}m_2^{\sigma}} \left(\frac{m_1}{m_2}\right)^{it}.
\end{eqnarray*}

Hence we get
\begin{eqnarray}\label{J1-mid}
J_{31}(T;\sigma)&&=T\sum_{m=1}^\infty\frac{-a_{3,j}^{(1)}(m;\sigma,T)a_{3,j}(m;\sigma,T) }{m^{2\sigma}}+O(\mathcal{S}(T)),\\
\mathcal{S}(T)&& :=\sum_{\stackrel{m_1,m_2\in {\Bbb N}}{m_1\not= m_2}} \frac{|a_{3,j}^{(1)}(m_1;\sigma,T)a_{3,j}(m_2;\sigma,T)|}
 {m_1^{\sigma}m_2^{\sigma}|\log \frac{m_1}{m_2}|}. \nonumber
\end{eqnarray}

Similar to (\ref{main-term-cons}), we have
\begin{eqnarray}\label{J3-cons-mid}
   \sum_{m=1}^\infty\frac{-a_{3,j}^{(1)}(m;\sigma,T)a_{3,j}(m;\sigma,T) }{m^{2\sigma}}
=B_{3,j}(\sigma)+O(T^{2-2\sigma}(\log T)^{(3+j)^2+1}),
\end{eqnarray}
where
\begin{equation}\label{J3-cons}
B_{3,j}(\sigma):=\sum_{m=1}^\infty\frac{-a_{3,j}^{(1)}(m;\sigma)a_{3,j}(m;\sigma) }{m^{2\sigma}}.
\end{equation}

Now we bound $\mathcal{S}(T).$ Write
\begin{equation}\label{ST-division}
\mathcal{S}(T)=\mathcal{S}_1(T)+\mathcal{S}_2(T),
\end{equation}
where
\begin{eqnarray*}
&&\mathcal{S}_1(T)=\sum_{\stackrel{1\leq m_1\not= m_2 }{|\log \frac{m_1}{m_2}|\geq 1/2}} \frac{|a_{3,j}^{(1)}(m_1;\sigma,T)a_{3,j}(m_2;\sigma,T)|}
 {m_1^{\sigma}m_2^{\sigma}|\log \frac{m_1}{m_2} |},\\
 &&\mathcal{S}_2(T)=\sum_{\stackrel{1\leq m_1\not= m_2 }{|\log \frac{m_1}{m_2} |< 1/2}} \frac{|a_{3,j}^{(1)}(m_1;\sigma,T)a_{3,j}(m_2;\sigma,T)|}
 {m_1^{\sigma}m_2^{\sigma}|\log \frac{m_1}{m_2} |}.
\end{eqnarray*}

We have
\begin{eqnarray}\label{S1T-bound}
\mathcal{S}_1(T)&&\ll \left(\sum_{m_1=1}^\infty\frac{|a_{3,j}^{(1)}(m_1;\sigma,T)}{m_1^\sigma}\right)
\left(\sum_{m_2=1}^\infty\frac{|a_{3,j}(m_2;\sigma,T)}{m_2^\sigma}\right)\\
&&\ll \left(\sum_{n_1\leq T, \ell_1\geq 1} \frac{|\tau_3^{(1)}(n_1)\mu(\ell_1)|}{n_1^\sigma\ell_1^{2-\sigma}}\right)
\left(\sum_{n_2\leq T, \ell_2\geq 1} \frac{|\tau_3(n_1)\mu(\ell_2)|}{n_2^\sigma\ell_2^{2-\sigma}}\right)\nonumber\\
&&\ll T^{2-2\sigma}(\log T)^5,\nonumber
\end{eqnarray}
where in the last step we used the bound $\sum_{n\leq x}\tau_3(n)\ll x(\log x)^2$ and partial summation.

It remains to bound $\mathcal{S}_2(T).$ The condition
$|\log \frac{m_1}{m_2} |<1/2$ implies that $m_1\asymp m_2.$
By Lagrange mean value theorem there exists some $u\asymp m_1\asymp m_2$ such that
$$|\log \frac{m_1}{m_2} |=u^{-1}|  m_1- m_2|\asymp (m_1m_2)^{-1/2}|m_1-m_2|.$$
Hence by the formula (5.10) of Ivi\'c \cite{Iv} we get
\begin{eqnarray}\label{S2T-bound}
\mathcal{S}_2(T)&&\ll \sum_{\stackrel{m_1\asymp m_2 \leq T }{ m_1\not= m_2}} \frac{|a_{3,j}^{(1)}(m_1;\sigma,T)a_{3,j}(m_2;\sigma,T)|}
 {m_1^{\sigma-1/2}m_2^{\sigma-1/2}|m_1-m_2|}\\
 &&\ll \left(\sum_{m_1=1}^\infty\frac{(a_{3,j}^{(1)}(m_1;\sigma,T))^2 }
 {m_1^{2\sigma-1} }\right)^{1/2}
 \left(\sum_{m_2=1}^\infty\frac{(a_{3,j}(m_2;\sigma,T))^2 }
 {m_2^{2\sigma-1} }\right)^{1/2}\nonumber\\
 &&\ll T^{2-2\sigma}(\log T)^{(3+j)^2}\nonumber
\end{eqnarray}
with the help of the bound (\ref{second-series}) and  the estimate
$$\sum_{m_2=1}^\infty\frac{(a_{3,j}(m_2;\sigma,T))^2 }
 {m_2^{2\sigma-1} }\ll T^{2-2\sigma}(\log T)^{(3+j)^{2}-1},$$
which can be proved by the same way as (\ref{second-series}).

Combining (\ref{J1-mid})-(\ref{S2T-bound}) we get
\begin{equation}\label{J1-zuihou}
J_{31}(T;\sigma)=B_{3,j}(\sigma)T+O(T^{2-2\sigma}(\log T)^{(3+j)^{2}+1}).
\end{equation}

We now consider $J_{32}(T;\sigma).$ By Cauchy's inequality and Lemma 2.4 we have
\begin{eqnarray} \label{S32-bound-mid}
 J_{32}(T,\sigma)&& \ll \left(\int_1^T|\zeta(\sigma+it)|^6dt\right)^{1/2}\\
&&\ \ \   \ \times
  \left(\int_1^T \left|\zeta^2(\sigma-it)\zeta^{\prime}(\sigma-it)
 +\sum_{m\leq T}  \frac{\tau_3^{(1)}(m)}{m^{\sigma-it}} \right|  ^2dt\right)^{1/2} \nonumber\\
&&\ll T^{\frac 12+\varepsilon} \left(\int_1^T \left|\zeta^2(\sigma-it)\zeta^{\prime}(\sigma-it)
 +\sum_{m\leq T}  \frac{\tau_3^{(1)}(m)}{m^{\sigma-it}} \right|^2dt\right)^{1/2}.\nonumber
\end{eqnarray}

We only need to bound
the integral
$$ \int_1^T \left|\zeta^2(\sigma-it)\zeta^{\prime}(\sigma-it)
 +\sum_{m\leq T}  \frac{\tau_3^{(1)}(m)}{m^{\sigma-it}} \right|^2dt
$$
by Lemma 2.12. This time we take $q=2, \alpha=7/12+\delta, \beta=1+\delta, \gamma=\sigma,$ where $0<\delta<1/2$ is a  positive constant which is arbitrarily
small and
$$F(s)= \zeta^2(s)\zeta^{\prime}(s)
 +\sum_{m\leq T}  \frac{\tau_3^{(1)}(m)}{m^{s}}.$$
By Lemma 2.4 and Lemma 2.11 we have
\begin{equation}\label{alpha-1}
\int_1^T|F(\alpha-it)|^2dt\ll T^{1+\delta}.
\end{equation}

From (\ref{diri-series}) we have
$$F(\beta-it)=-\sum_{n>T}\frac{\tau_3^{(1)}(n)}{n^{1+\delta-it}}.$$
Thus
$$F^2(\beta-it)=\sum_{n>T^2}\frac{b(n,T)}{n^{1+\delta-it}},$$
where $|b(n,T)|\leq \tau_6(n)\log^2 n.$ So by Cauchy's inequality and Lemma 2.11
we get
\begin{eqnarray}\label{beta-1}
\int_1^T|F(\beta-it)|^2dt&&\ll T^{1/2}
\int_1^T \left|\sum_{n>T^2}\frac{b(n,T)}{n^{1+\delta-it}}\right|^2dt\ll T^{1/2}.
\end{eqnarray}

From (\ref{alpha-1}), (\ref{beta-1}) and Lemma 2.12 we get
\begin{eqnarray}\label{beta-2}
\int_1^T|F(\sigma-it)|^2dt \ll T^{\frac{17-12\sigma}{10}+\varepsilon}
\end{eqnarray}
if $\delta<\varepsilon/3.$

From (\ref{S32-bound-mid}) and (\ref{beta-2}) we get
\begin{equation}\label{J2-bound}
J_{32}(T;\sigma)\ll T^{(27-12\sigma)/20+\varepsilon}.
\end{equation}

Similarly we can prove
\begin{equation}\label{J3-bound}
J_{33}(T;\sigma)\ll T^{ (27-12\sigma)/20+\varepsilon}.
\end{equation}

From  (\ref{T3j-division}), (\ref{J1-zuihou}),
(\ref{J2-bound}) and (\ref{J3-bound}) we get
\begin{equation}\label{I3j-asymp}
I_{3,j}(T;\sigma)=B_3,j(\sigma)T
+O(  T^{(27-12\sigma)/20+\varepsilon}),\ \ \ \ \ (7/12<\sigma<1),
\end{equation}
where $B_{3,j}(\sigma)$ was defined in (\ref{J3-cons}).

\section{Proofs of Theorem 3 and corollaries}

In this section we prove Theorem 3 and two corollaries.

\subsection{A weighted mean value of $\Delta_{\mathcal{T}_{1,1}}(x)$  }\

 In this subsection we study the weighted mean value of $\Delta_{\mathcal{T}_{1,1}}(x), $
which was defined in (\ref{error-j}).

Suppose $z>0.$
From (3.1) and  Perron's formula we have
\begin{eqnarray}\label{error}
 \Delta_{\mathcal{T}_{1,1}}(z)=\frac{1}{2\pi i}\lim_{T\rightarrow \infty}\int_{\lambda-iT}^{\lambda+iT}
-\frac{\zeta^2(s-1)\zeta^{\prime}(s-1)}{\zeta(s)}z^ss^{-1}ds,
\end{eqnarray}
where $\lambda<2$ and close to $2.$ From Lemma 2.4 we see
that (\ref{error}) holds for any $\lambda>4/3.$ Replacing in (\ref{error}) $z$ by $1/z$, taking $4/3<\lambda<3/2,$ and using Parseval's identity (see formula (A5)  in the appendix of  Ivi'c  \cite{Iv} )
we get
\begin{eqnarray}\label{relation}
\ \ \ \ \  \frac{1}{2\pi}\int_{-\infty}^\infty
\frac{|\zeta^2(\lambda-1+it)\zeta^{\prime}(\lambda-1+it)|^2}{|\zeta(\lambda+it)|^2|\lambda+it|^2}dt
= \int_0^\infty \frac{( \Delta_{\mathcal{T}_{1,1}}(z))^2}{z^{2 \lambda+1}}dz.
\end{eqnarray}

Write
\begin{eqnarray}\label{transss}
\int_0^\infty \frac{( \Delta_{\mathcal{T}_{1,1}}(z))^2}{z^{2 \lambda+1}}dz &&=\int_1^\infty \frac{( \Delta_{\mathcal{T}_{1,1}}(z))^2}{z^{2 \lambda+1}}dz +O(1)\\
&&=\int_1^\infty \frac{ ( \Delta_{\mathcal{T}_{1,1}}(z))^2}{z^{2( \lambda-4/3)+11/3}}dz+O(1)\nonumber\\
&&=\int_1^\infty e^{-\xi x}dA(x)+O(1)\nonumber
\end{eqnarray}
by the changes of variables $x=2\log z$ and $\xi=\lambda-4/3,$
where
\begin{equation}\label{Ax-def}
A(x):=\int_1^{e^{x/2}}\frac{ ( \Delta_{T_{1,1}}(u))^2 }{u^{11/3}}du.\end{equation}

 Define
\begin{eqnarray} \label{H-sigma-def}
\ \ \ H(\xi):=\frac{1}{2\pi}
\int_{-\infty}^\infty
\frac{|\zeta^2(\xi+1/3+it)\zeta^{\prime}(\xi+1/3+it)|^2}
{|\zeta(\xi +4/3+it)|^2|\xi+4/3+it|^2}dt,\ (0<\xi<1/12).
\end{eqnarray}
 Then
\begin{eqnarray}\label{H-sigma}
H(\xi)&&=\frac{1}{\pi}
\int_1^\infty
\frac{|\zeta^2(\xi+1/3+it)\zeta^{\prime}(\xi+1/3+it)|^2}
{ |\zeta(\xi+4/3+it)|^2 |\xi+4/3+it|^2}dt
+O(1)\\
&&=\frac{1}{4\pi^3}
\int_1^\infty
\frac{|\zeta^2(\xi+1/3+it)\zeta^{\prime}(\xi+1/3+it)|^2}
{|\zeta(\xi+4/3+it)|^2\times (t/2\pi)^2}dt
+O(1),\nonumber
\end{eqnarray}
by noting that
$$\int_1^\infty
 \frac{|\zeta^2(\xi+1/3+it)\zeta^{\prime}(\xi+1/3+it)|^2}
  {t^3}dt\ll 1,$$
 which follows easily from Lemma 2.3.

 From Lemma 2.1 and Lemma 2.2 we have ($s=\sigma+it, 0<\sigma<1/2,\ \ t\geq 1$)
\begin{eqnarray*}
-\zeta^2(s)\zeta^{\prime}(s)&&=\chi^3(s)\zeta^2(1-s)\zeta^{\prime}(1-s)-
\chi^2(s)\chi^{\prime}(s)\zeta^3(1-s) \\
&&=\chi^3(s)\left(\zeta^2(1-s)\zeta^{\prime}(1-s)-
 \frac{\chi^{\prime}(s)}{\chi(s)}\zeta^3(1-s)\right)\\
 &&=\chi^3(s)\left(\zeta^2(1-s)\zeta^{\prime}(1-s)
 + \zeta^3(1-s)\log\frac{t}{2\pi}+O\left(\frac{|\zeta^3(1-s)|}{t}\right)\right).
\end{eqnarray*}
Thus we have for $s= \xi+1/3+it\ \ (0<\xi<1/12,\ t\geq 1)$ that
\begin{eqnarray}\label{function-division}
&&\ \ \ \ \ \frac{|\zeta^2(s)\zeta^{\prime}(s)|^2}{|\zeta(1+s)|^2} \ (t/2\pi)^{-2}\\
&&=\frac{|\chi(s)|^6}{|\zeta(1+s)|^2}\left(\zeta^2(1-s)\zeta^{\prime}(1-s)
 + \zeta^3(1-s)\log\frac{t}{2\pi}+O\left(\frac{|\zeta^3(1-s)|}{t}\right)\right)
 \nonumber\\&&
 \ \ \ \ \ \ \ \ \ \ \ \ \ \ \times   \overline{\left(\zeta^2(1-s)\zeta^{\prime}(1-s)
 + \zeta^3(1-s)\log\frac{t}{2\pi}+O\left(\frac{|\zeta^3(1-s)|}{t}\right)\right)}
 \nonumber\\
 &&=   \frac{|\zeta^3(2/3-\xi+it)|^2}{|\zeta(\xi+4/3+it)|^2}\times \left(\frac{t}{2\pi}\right)^{-1-6\xi}\log^2\frac{t}{2\pi}\nonumber\\
 &&\ \ \ \  +
   \frac{|\zeta^2(2/3-\xi+it)  \zeta^{\prime}(2/3-\xi+it)|^2}
   {|\zeta(\xi+4/3+it)|^2} \times \left(\frac{t}{2\pi}\right)^{-1-6\xi} \nonumber\\
 &&\ \ \ \ + \frac{\zeta^2(2/3-\xi-it)\zeta^{\prime}(2/3-\xi-it)}
 {\zeta(\xi+4/3-it)}\frac{
 \zeta^3(2/3-\xi+it)}{\zeta(\xi+4/3+it)}\times  \left(\frac{t}{2\pi}\right)^{-1-6\xi}\log\frac{t}{2\pi}\nonumber\\
 &&\ \ \ \ +  \frac{\zeta^2(2/3-\xi+it)\zeta^{\prime}(2/3-\xi+it)}{\zeta(\xi+4/3+it)}
 \frac{\zeta^3(2/3-\xi-it)}{\zeta(\xi+4/3-it)}
 \left(\frac{t}{2\pi}\right)^{-1-6\xi}\log\frac{t}{2\pi}\nonumber\\
 &&\ \ \ +O\left(\frac{|\zeta(2/3-\xi-it)|^6\log t}{t^{2+6\xi}}
 +\frac{|\zeta(2/3-\xi+it)|^5 |\zeta^{\prime}(2/3-\xi+it)|}{t^{2+6\xi}} \right)\nonumber\\
 &&=\sum_{j=1}^5 h_j(t;\xi),\nonumber
\end{eqnarray}
say.

From (\ref{H-sigma}) and {\ref{function-division}) we have
\begin{eqnarray}
H(\xi)=\frac{1}{4\pi^3}\sum_{j=1}^5 \mathcal{I}_j(\xi),
\end{eqnarray}
where
$$\mathcal{I}_j(\xi):=\int_1^\infty h_j(t;\xi)dt,\ \ j=1, 2, 3, 4, 5. $$

Let $\alpha>0$ be a fixed real number and $k\geq 0$ be a fixed integer. Define
$$J_k(\alpha):=\int_1^\infty\frac{\log^k t}{t^{1+\alpha}}dt.$$
It is easy to see that $J_0(\alpha)=1/\alpha.$ By partial integration it is easy to see that
\begin{equation}\label{J-k}
J_k(\alpha)=\frac{k!}{\alpha^{k+1}}, \ \ k\geq 0.
\end{equation}

We first evaluate $\mathcal{I}_1(\xi).$
Let $$E_{I_{2,1}}(T;\sigma):=I_{2,1}(T;\sigma)-B_{2,1}(\sigma)T
\ \ \ (\sigma=2/3-\xi),$$
where $I_{2,1}(T;\sigma)$ was defined in (\ref{integrals}), $B_{2,1}(\sigma)$ was defined in (\ref{I2j-cons}).
From (\ref{222-final}) we have
\begin{equation}\label{I2-error}
E_{I_{2,1}}(T;\sigma)\ll T^{1-(12\sigma-7)/10+\varepsilon}=T^{1-(1-12\xi)/10+\varepsilon}.
\end{equation}
By partial integration, we can write
\begin{eqnarray}\label{I1+I2}
\mathcal{I}_1(\xi)&&=\int_1^\infty \left(\frac{t}{2\pi}\right)^{-1-6\xi}\log^2\frac{t}{2\pi}
d\int_1^t\frac{|\zeta(2/3-\xi+iu)|^6}{|\zeta(\xi+4/3+iu)|^2}du\\
&&=\int_1+\int_2,\nonumber
\end{eqnarray}
where
\begin{eqnarray*}
&&\int_1:=B_{2,1}\left(\frac 23-\xi\right)\int_1^\infty \left(\frac{t}{2\pi}\right)^{-1-6\xi}\log^2\frac{t}{2\pi}
dt,\\
&&\int_2:=\int_1^\infty \left(\frac{t}{2\pi}\right)^{-1-6\xi}\log^2\frac{t}{2\pi}
dE_{I_{2,1}}(t;2/3-\xi).
\end{eqnarray*}

By the definition of $J_k(\alpha)$ and (\ref{J-k}) we get
\begin{eqnarray}\label{I1-asy}
\ \ \ \int_1&&=B_{2,1}\left(\frac 23-\xi\right)\int_1^\infty \left(\frac{t}{2\pi}\right)^{-1-6\xi}\log^2\frac{t}{2\pi}dt \\
&&= B_{2,1}\left(\frac 23-\xi\right)(2\pi)^{1+6\xi} \int_1^\infty t^{-1-6\xi}
(\log t-\log 2\pi)^2dt\nonumber\\
&&= B_{2,1}\left(\frac 23-\xi\right)(2\pi)^{1+6\xi} \left(J_2(6\xi)-2 J_1(6\xi)\log 2\pi
+J_0(6\xi)\log^2 2\pi\right)\nonumber\\
&&= B_{2,1}\left(\frac 23-\xi\right)(2\pi)^{1+6\xi} \left(\frac{1}{108\xi^3} -
\frac{\log 2\pi}{18\xi^2} +\frac{\log^2 2\pi}{6\xi}
 \right)\nonumber\\
 &&=\frac{B_{2,1}(2/3)\pi}{54}\times\frac{1}{\xi^3}
 +O\left(\frac{1}{\xi^2}\right)\ \ (\xi \rightarrow 0)\nonumber
\end{eqnarray}
by noting that $B_{2,1}(\sigma)$ is a continuous function on the interval $(7/12,1)$ and that
$$B_{2,1}(\sigma)=B_{2,1}(2/3)+O(|\sigma-2/3|)\ \ \ \ (\sigma\rightarrow 2/3).$$

By (\ref{I2-error}) and partial integration it is easy to check that
\begin{eqnarray}\label{I2-upper-bound}
\int_2\ll 1,\ \ \ 0<\xi<1/12.
\end{eqnarray}

From (\ref{I1+I2})-(\ref{I2-upper-bound}) we get  that
\begin{eqnarray}\label{I1-asssp}
\mathcal{I}_1(\xi)  =\frac{B_{2,1}(2/3)\pi}{54}\times\frac{1}{\xi^3}+O\left(\frac{1}{\xi^2}\right),\ \
\xi\rightarrow 0^+.
\end{eqnarray}

Similarly, from (\ref{111-final}), (\ref{222-final}), (\ref{J-k}) and partial summation  we get the following estimates:
\begin{eqnarray}
\mathcal{I}_2(\xi)  =\frac{B_{1,1}(2/3)\pi}{3}\times\frac{1}{\xi}+O(1), \ \ \xi\rightarrow 0^+.
\end{eqnarray}

\begin{eqnarray}
\mathcal{I}_3(\xi) = \mathcal{I}_4(\xi) =\frac{B_{3,1}(2/3)\pi}{18}\times\frac{1}{\xi^2}+O\left(\frac{1}{\xi}\right),\ \ \xi\rightarrow 0^+.
\end{eqnarray}

 From Lemma 2.4 and partial integration we get easily that
 \begin{eqnarray}\label{I5-bun}
\mathcal{I}_5(\xi) \ll 1, \ \ -1/12<\xi<1/12.
\end{eqnarray}

Combining (\ref{I1-asssp})-(\ref{I5-bun}) we get
\begin{eqnarray}\label{Hxi-asp}
H(\xi)  =\frac{B_{2,1}(2/3)}{216\pi^2}\times\frac{1}{\xi^3}+O\left(\frac{1}{\xi^2}\right),\ \ \xi\rightarrow 0^+.
\end{eqnarray}

From (\ref{relation}), (\ref{transss}) and   (\ref{H-sigma-def}) we have
\begin{equation}\label{relation2}
H(\xi)=\int_1^\infty e^{-\xi x}dA(x)+O(1)
\end{equation}
 where $A(x)$ was defined in (\ref{Ax-def}).

From (\ref{Ax-def}),  (\ref{Hxi-asp}), (\ref{relation2}) and Lemma 2.13 by taking
$C=B_{2,1}(2/3)/216\pi^2,\ \psi(x)=x,\ \omega=3$ we get
$$A(x)=\int_1^{e^{x/2}}\frac{ ( \Delta_{\mathcal{T}_{1,1}}(u))^2 }{u^{11/3}}du
=\left(\frac{B_{2,1}(2/3)}{216\pi^2}+O\left(\frac{1}{\log x}\right)\right)\frac{x^3}{\Gamma(4)},
$$
 which implies that
 \begin{eqnarray}\label{delta-t1}
 \int_1^X\frac{ ( \Delta_{\mathcal{T}_{1,1}}(u))^2 }{u^{11/3}}du
 =\left(\frac{B_{2,1}(2/3)}{162\pi^2}+O\left(\frac{1}{\log\log X}\right)\right) \log^3 X.
 \end{eqnarray}

Similar to (\ref{delta-t1}), we can prove that
\begin{eqnarray}\label{delta-t2}
 \int_1^X\frac{ ( \Delta_{\mathcal{T}_{2,1}}(u))^2 }{u^{11/3}}du
 =\left(\frac{B_{2,1}(2/3)}{3\pi^2}+O\left(\frac{1}{\log\log X}\right)\right) \log X.
 \end{eqnarray}

\subsection{The mean-square of $\Delta_{\mathcal{T}_{3,1}}(x)$}\

In this subsection, we study the mean-square of $\Delta_{\mathcal{T}_{3,1}}(x),$
which was defined by (\ref{error-j}) for $i=3.$

Similar to (\ref{relation}), we have the expression
\begin{eqnarray}\label{relation-3333}
\ \ \ \ \  \frac{1}{2\pi}\int_{-\infty}^\infty
\frac{|\zeta(\lambda-1+it)R(\lambda+it)|^2}{|\zeta(\lambda+it)|^2|\lambda+it|^2}dt
= \int_0^\infty \frac{( \Delta_{\mathcal{T}_{3,1}}(z))^2}{z^{2 \lambda+1}}dz,
\end{eqnarray}
where $R(s)$ was defined in (\ref{Rs-def}).

On the line $\Re s=1+\varepsilon+it,$ by Lemma 2.3 and (\ref{Rs-bound}) we have
$$\zeta(s-1)\ll (|t|+1)^{1/2-2\varepsilon/3},\ \ \ R(s)\ll (|t|+1)^{1/2-2\varepsilon}\log (|t|+2).$$
From (\ref{relation-3333}) we see that
\begin{eqnarray}\label{relation-33}
 \int_0^\infty \frac{( \Delta_{\mathcal{T}_{3,1}}(z))^2}{z^{2 \lambda+1}}dz
 \ll 1, \ \ \ \lambda\geq 1+\varepsilon,
\end{eqnarray}
which implies that
\begin{eqnarray*}\label{delta-3-mean}
 \int_1^X ( \Delta_{\mathcal{T}_{3,1}}(z))^2 dz\ll X^{3+\varepsilon}.
\end{eqnarray*}

\subsection{Proof of Theorem 3}\

We first consider the case $S^{(1)}(x).$ From (\ref{S1-error}), we have
\begin{eqnarray}\label{mix}
&&\ \ \ \ \ \ \frac{(E_{S^{(1)}}(x))^2}{x^{11/3}}\\&&=
16\frac{(\Delta_{\mathcal{T}_{1,1}}(x))^2}{x^{11/3}}+(32\gamma-16)
\frac{\Delta_{\mathcal{T}_{1,1}}(x)\Delta_{\mathcal{T}_{2,1}}(x) }{x^{11/3}}
  +(4\gamma-2)^2\frac{(\Delta_{\mathcal{T}_{2,1}}(x))^2}{x^{11/3}}\nonumber\\
  &&\ \ +O\left(\frac{|\Delta_{\mathcal{T}_{1,1}}(x)\Delta_{\mathcal{T}_{3,1} }(x)| }{x^{11/3}}\right)+O\left(\frac{|\Delta_{\mathcal{T}_{2,1}}(x)\Delta_{\mathcal{T}_{3,1} }(x)| }{x^{11/3}}\right)+\frac{(\Delta_{\mathcal{T}_{3,1}}(x))^2}{x^{11/3}}\nonumber\\
  &&\ \ +O\left(\frac{|\Delta_{\mathcal{T}_{1,1}}(x)|}{x^{8/3-\varepsilon}}\right)
  +O\left(\frac{|\Delta_{\mathcal{T}_{2,1}}(x)|}{x^{8/3-\varepsilon}}\right)
  +O\left(\frac{|\Delta_{\mathcal{T}_{3,1}}(x)|}{x^{8/3-\varepsilon}}\right).\nonumber
\end{eqnarray}

From (\ref{delta-t1}), (\ref{delta-t2}) and Cauchy's inequality we get
\begin{equation}\label{mixed}
\int_1^X \frac{\Delta_{\mathcal{T}_{1,1}}(x)\Delta_{\mathcal{T}_{2,1}}(x) }{x^{11/3}}dx\ll \log^2 X.
\end{equation}

From (\ref{delta-t1}), (\ref{delta-t2}),  (\ref{relation-33}) and some easy calculations we have the following estimates
\begin{eqnarray}\label{mixed-mix}
&&\int_1^X \frac{|\Delta_{\mathcal{T}_{1,1}}(x)\Delta_{\mathcal{T}_{3,1} }(x)| }{x^{11/3}}   dx\ll 1, \ \ \
 \int_1^X \frac{|\Delta_{\mathcal{T}_{2,1}}(x)\Delta_{\mathcal{T}_{3,1} }(x)| }{x^{11/3}}   dx\ll 1, \\
&&\int_1^X \frac{|\Delta_{\mathcal{T}_{1,1}}(x)|}{x^{8/3-\varepsilon}}dx\ll 1,\ \ \ \
\int_1^X \frac{|\Delta_{\mathcal{T}_{2,1}}(x)|}{x^{8/3-\varepsilon}}dx\ll 1,\nonumber\\&&\int_1^X \frac{|\Delta_{\mathcal{T}_{3,1}}(x)|}{x^{8/3-\varepsilon}}dx\ll 1.\nonumber
\end{eqnarray}

From  (\ref{delta-t1}), (\ref{delta-t2}) and  (\ref{relation-33})-(\ref{mixed-mix}) we get
\begin{equation*}
\int_1^X \frac{(E_{S^{(1)}}(x))^2}{x^{11/3}}dx=
\left(\frac{8B_{2,1}(2/3)}{81\pi^2}+O\left(\frac{1}{\log\log X}\right)\right) \log^3 X.\end{equation*}

This completes the proof of Theorem 3 for $f=S^{(1)}.$
By the same approach, for $f=C^{(1)}$ we have
\begin{equation*}
\int_1^X \frac{(E_{C^{(1)}}(x))^2}{x^{11/3}}dx=
\left(\frac{8B_{2,2}(2/3)}{81\pi^2}+O\left(\frac{1}{\log\log X}\right)\right) \log^3 X.\end{equation*}

For $S^{(2)}(x),$ we have the formula (\ref{S2-error-expr}). From Proposition 4.2 we see that the average order of $E_U(x)$ is $x^{5/4}.$ So repeating the proof of the case
$S^{(1)}(x),$ we get
\begin{equation*}
\int_1^X \frac{(E_{S^{(2)}}(x))^2}{x^{11/3}}dx=
\left(\frac{8B_{2,1}(2/3)}{81\pi^2}+O\left(\frac{1}{\log\log X}\right)\right) \log^3 X.\end{equation*}

Similarly, we have
\begin{equation*}
\int_1^X \frac{(E_{C^{(2)}}(x))^2}{x^{11/3}}dx=
\left(\frac{8B_{2,2}(2/3)}{81\pi^2}+O\left(\frac{1}{\log\log X}\right)\right) \log^3 X.\end{equation*}

This completes the proof of Theorem 3.

\subsection{Proof of Corollary 1}\

The  $\Omega$-result  (\ref{1.7}) in Corollary 1 is a direct consequence of Theorem 3. So we only prove the estimate (\ref{1.6}). Suppose $f(x)\in \mathcal{F}(x).$
Let
$$\mathcal{H}_f(x):=\int_1^x\frac{E_f^2(t)}{t^{11/3}}dt.$$

From Theorem 3 and partial integration we have
\begin{eqnarray*}
&&\ \ \ \int_1^XE_f^2(x)dx=\int_1^X\frac{E_f^2(x)}{x^{11/3}}x^{11/3}dx=
\int_1^Xx^{11/3}d\mathcal{H}_f(x)\\
&&=X^{11/3}\mathcal{H}_f(X)-\frac{11}{3}\int_1^Xx^{8/3}\mathcal{H}_f(x)dx\nonumber\\
&&=c_fX^{11/3}\log^3 X-\frac{11c_f}{3}\int_1^Xx^{8/3}\log^3 xdx
+O\left(\frac{X^{11/3}\log^3 X}{\log\log X}\right)\nonumber\\
&&=3c_f\int_1^X x^{8/3}\log^2 xdx+O\left(\frac{X^{11/3}\log^3 X}{\log\log X}\right)\nonumber\\
&&\ll \frac{X^{11/3}\log^3 X}{\log\log X}.\nonumber
\end{eqnarray*}

\subsection{Proof of Corollary 2}\

In this subsection we prove Corollary 2.

\subsubsection{Proof of the first assertion in Corollary 2}\

Suppose $f(x)\in \mathcal{F}(x).$   We shall use reduction to absurdity.

From Theorem 3 there exists a positive constant $D>0$  such that
\begin{equation}\label{TH3-inequality}
\left|\int_1^X\frac{E_f^2(x)}{x^{11/3}}dx-c_f\log^3 X\right|\leq \frac{D\log^3 X}{\log\log X},\ \ X\geq 100.
\end{equation}

Suppose $T\geq 100$ is any large parameter. From (\ref{TH3-inequality}) it is easy to check that
\begin{equation}\label{y-range}
\int_T^{yT}\frac{E_f^2(x)}{x^{11/3}}dx\gg \log^2 T\times \log y
\end{equation}
if $\log y\geq \frac{D\log T}{\log\log T}.$
Now take $y=y(T)=e^{\frac{D\log^3 T}{\log\log T}}.$ Then $y\ll T^\varepsilon.$ Suppose
$E_f(x)$ doesn't have sign changes in $(T, yT]$. Then we have $E_f(x)>0$ or
$E_f(x)<0$ in this interval. Without loss of generality,
 suppose $E_f(x)>0.$ Then from Theorem 1 and  Theorem 2 we get
\begin{eqnarray*}
\int_T^{yT}\frac{(E_{f}(x))^2}{x^{11/3}}&&\ll T^{\frac{139}{96}-\frac{11}{3}+4\varepsilon}
\int_T^{yT} E_{f}(x)dx\\
&& \ll T^{\frac{139}{96}-\frac{11}{3}+\frac{13}{6}+6\varepsilon}\ll
T^{-\frac{5}{96}+6\varepsilon},
\end{eqnarray*}
which is contrary to (\ref{y-range}). This means that {\it $E_f(x)>0$ in $(T, yT]$} is not true. Similarly {\it $E_f(x)<0$ in $(T, yT]$} is not true. So $E_f(x)$ changes its sign in  $(T, yT]$ at least once.
 Since $T$ is arbitrary, we see that $E_f(x)$ must have infinitely many sign changes in $(1, \infty)$.

\subsubsection{Proof of the second assertion in Corollary 2}\

 Let $T_0=\max(100, e^{e^{2D}}),$ where $D>0$ is the constant defined in (\ref{TH3-inequality}). We define a sequence $\{T_n\}(n\geq 1)$ by the recurrence formula
\begin{equation}\label{recurrence}
T_{n+1}=T_n\times e^{\frac{D\log T_n}{\log\log T_n}}, \ \ n\geq 0.
\end{equation}
From the proof in Section 8.4.1 we see that $E_f(x)$
has at least one sign change in the interval $(T_n, T_{n+1}]\ (n\geq 0).$
It is easy to see that the sequence $\{T_n\}\ (n\geq 0)$ is increasing, tends to infinity when $n$ tends to infinity and that for each $n,$ the estimate
\begin{equation}\label{neibour}
T_n<T_{n+1}\ll T_n^{1+\varepsilon}
\end{equation}
holds.

Let $X>100$ be a large parameter.
Then there exists a unique integer $\mathcal{L}\geq 1$    such that
$T_{\mathcal{L}}<X\leq T_{\mathcal{L}+1}.$ Obviously we have $N_f(X)\geq \mathcal{L}.$
From (\ref{recurrence}) we see that
\begin{equation}\label{quotion}
\frac{\log T_{n+1}}{\log T_n }= 1+\frac{D }{\log\log T_n}, \ \ (n\geq 0).
\end{equation}

Define an integer $\mathcal{L}_0$ such that
\begin{equation}\label{upper-lower}
T_{\mathcal{L}_0}\leq e^{(\log X)^{1/2}}<T_{\mathcal{L}_0+1}.
\end{equation}

From (\ref{quotion}) we get
$$\frac{\log T_{\mathcal{L}+1}}{\log T_{\mathcal{L}_0} }
=\prod_{n=\mathcal{L}_0}^{\mathcal{L}}\left(1+\frac{D }{\log\log T_n}\right),$$
which is equivalent to
\begin{eqnarray}\label{left-side}
\log\frac{\log T_{\mathcal{L}+1}}{\log T_{\mathcal{L}_0} } =
 \sum_{n=\mathcal{L}_0}^{\mathcal{L}}\log\left(1+\frac{D }{\log\log T_n}\right).
\end{eqnarray}

From $T_{\mathcal{L}}<X\leq T_{\mathcal{L}+1}$ and (\ref{upper-lower}) we have
$ \frac{\log T_{\mathcal{L}+1}}{\log T_{\mathcal{L}_0} }\geq (\log X)^{1/2}.$
 Thus we get
 \begin{eqnarray}
\log\frac{\log T_{\mathcal{L}+1}}{\log T_{\mathcal{L}_0} } \geq \frac 12\log\log X.
\end{eqnarray}

 From $\log(1+t)\ll t\ (0<t<1/2)$ and the increasing property of $T_n$ we get
 \begin{eqnarray}
 \sum_{n=\mathcal{L}_0}^{\mathcal{L}}\log\left(1+\frac{D }{\log\log T_n}\right)
 \ll \sum_{n=\mathcal{L}_0}^{\mathcal{L}}\frac{1 }{\log\log T_n}
 \ll \frac{ \mathcal{L} }{\log\log T_{\mathcal{L}_0}}.
 \end{eqnarray}

 From (\ref{neibour}) and (\ref{upper-lower}) it is easy to check that
\begin{equation}\label{nnnnnnnn}
\log\log T_{\mathcal{L}_0}\asymp \log\log X.
\end{equation}

Gathering  (\ref{left-side})-(\ref{nnnnnnnn}) we get that  $N_f(X)\geq \mathcal{L}\gg (\log\log X)^2.$

\subsubsection{Proof of the third assertion in Corollary 2}\

Let $M\geq 100$ be any large real number. Suppose $M<T\leq 2M$ and $y=o(M)$ is a
parameter to be determined later. If  the asymptotic formula (\ref{con-asy-2}) were
true, then it is easy to see that
\begin{equation}\label{lower-bound-bound}
\int_T^{T+y}E_f^2(x)dx\gg yT^{8/3}\log^2 T
\end{equation}
if $y\geq  M^{1-\eta}$.

Take $y=\max(M^{1-\eta}, M^{91/96+2\varepsilon}).$ We shall show that
$E_f(x)$ has at least one sign change in the interval $(T, T+y].$
If $E_f(x)>0$ for all $x\in (T, T+y],$ then from Theorem 1 and Theorem 2 again we have
\begin{eqnarray*}
\int_T^{T+y}E_f^2(x)dx\ll T^{\frac{139}{96}+\varepsilon}\int_T^{T+y}E_f(x)dx
\ll T^{\frac{139}{96}+\frac{13}{6}+2\varepsilon}=T^{\frac{347}{96}+2\varepsilon},
\end{eqnarray*}
which is contrary to (\ref{lower-bound-bound}) by noting the choice of $y$. Hence we see that $E_f(x)$   has at least one sign change in $(T, T+y].$ Totally in the interval $(M, 2M]$, $E_f(x)$ has at least $M/y=\min(M^\eta, M^{5/96-2\varepsilon})$ sign changes. This estimate gives $N_f(X)\gg \min(X^\eta, X^{5/96-2\varepsilon}).$


\bigskip





\vskip 10mm

\begin{thebibliography}{100}



\bibitem{Cr}H. Cram\'er, \"Uber zwei S\"atze von Herrn G. H. Hardy,
Math. Z. {\bf 15}(1922), 201-210.


\bibitem{EST}D. Essouabri, C. Salinas Zavala, L. T\'{o}th,
Mean values of multivariable multiplicative functions and
applications to the average number of cyclic subgroups and
multivariable averages associated with the LCM function,
Journal of Number Theory {\bf 236}(2022), 404-442.



\bibitem{Fu}J. Furuya, On the summatory function of a product of an arithmetical function and its relevant error term, Ann. Sci. Math.
     Qu\'{e}bec {\bf 31}(2007), no.2, 165-185.

\bibitem{Fu1}J. Furuya,  On the average orders of the error term in the
Dirichlet divisor problem, J. of Number Theory. {\bf 115}(2005), 1-26.

\bibitem{HB}D. R. Heath-Brown, The distribution and moments of the error term in the Dirichlet divisor problem, Acta Arith. {\bf 60}(1992), 389--415.


\bibitem{Ko}G. Kolesnik, On the estimation of multiple esponential sums,
{\it Recent progress in analytic number theory}, Vol. 1(Durham, 1979), pp. 231-246,
Academic Press, London-New York, 1981.



\bibitem{Iv}A. Ivi\'c, The Riemann Zeta-Function. Theory and Applications, Wiley, New York, 1985.

\bibitem{Iv2}A. Ivi\'c, The theory of Hardy's Z-Function. Cambridge University Press.
2013.

\bibitem{Iv3}A. Ivi\'c, Large values of the error term in the divisor problem,
Invent. Math. {\bf 71}(1983), 513-520.

\bibitem{IZ}A. Ivi\'c and W. Zhai, Higher moments of the error term
in the divisor problem, Mathematical Notes, Vol. {\bf 88}, No. 3, (2010),
 pp. 338-346.

\bibitem{LT}Y. K. Lau and K.-M. Tsang, On the mean square formula of the error term in the Dirichlet divisor problem, Mathematical Proceedings of the Cambridge Philosophical Society, {\bf 146}(2009), 277-287.

\bibitem{NT} W. G. Nowak and L. T\'{o}th,  On the average number of subgroups of the group ${\Bbb Z}_{m}\times {\Bbb Z}_{n}$,
  Int. J. Number Theory {\bf 10} (2014), 363-374.

 \bibitem{KV} A. A. Karatsuba and S. M. Voronin, The Riemann zeta-function, Walter de Gruyter, Berlin--New York, 1992.

\bibitem{Ra}K. Ramachandra, On the mean-value and omega-theorems for the Riemann zeta-function, LN's {\bf 85}, Tata Inst. of Fundamental Research(distr. by Springer Verlag, Berlin etc.), Bombay, 1995.

\bibitem{Te} G. Tenenbaum, Introduction to analytic and probabilistic number theory,
Third edition, Graduate studies in Mathematics, {\it American Mathematical Society, RL}, 2015.


\bibitem{T}K. C. Tong, On divisor problem III, Acta Math. Sinica {\bf 6}
(1956), 515-541.

\bibitem{To}L. T\'{o}th,  On the number of cyclic subgroups of a finite Abelian group. Bull. Math. Soc. Sci. Math. Roumanie. vol. {\bf 55}(2012), 423-428.

\bibitem{TZ1}L. T\'{o}th and W. Zhai, On the error term concerning the number of subgroups of the groups ${\Bbb Z}_m\times {\Bbb Z}_n$ with $m,n\leq x,$ Acta
    Arith. {\bf 183} (2018), 285-299.

\bibitem{TZ2}L. T\'{o}th and W. Zhai, On the average number of cyclic subgroups  of the groups ${\Bbb Z}_{n_1}\times {\Bbb Z}_{n_2}\times {\Bbb Z}_{n_3}$ with $n_1, n_2, n_3\leq x,$ {\it Res. Number Theory}(2020) 6:12.

\bibitem{Ts} K.-M. Tsang, Higher-power moments of $\Delta(x), E(t)$ and $P(x)$,
Proc. London Math. Soc.(3) {\bf 65}(1992), 65-84.

\bibitem{Ts2}K.-M. Tsang, Recent progress on the Dirichlet divisor problem and the mean square of the Riemann zeta-function, Sci. China Math. {\bf 53}(2010), no. 9,
    2561-2572.

\bibitem{Vo}G. F. Vorono\"{\i}, Sur une fonction transcendante et ses applications \`{a} la sommation de quelques s\'{e}ries, Ann. \'{E}cole Normale (3) 21 (1904) 207-267, 459-533.

\bibitem{Zh1} W. Zhai, On higher-power moments of $\Delta(x)$ (II), Acta Arith.
{\bf 114}(2004), 35-54.


\end{thebibliography}
\end{document}